\documentclass[a4paper,12pt]{article}

\advance\textheight2.9cm        
\advance\topmargin-2.0cm
\advance\textwidth2.6cm
\advance\evensidemargin-1.6cm
\advance\oddsidemargin-1.6cm

\usepackage{amsfonts} 

\usepackage{url}

\usepackage{theorem}
\usepackage{amsfonts}
\usepackage{amssymb} 
\usepackage{bm}  
\usepackage{booktabs}
\def\comment#1{{}}	
\theoremstyle{change}
\newtheorem{thm}{Theorem.\nopagebreak}[section]

\newtheorem{prop}[thm]{Proposition.\nopagebreak}

\theorembodyfont{\rmfamily}

\newtheorem{remark}[thm]{Remark.\nopagebreak}

\def\Box{{\hbox{\raisebox{0.0em}{\rlap{$\sqcap$}}\kern0em%
			\raisebox{-0.0em}{$\sqcup$}}} } 

\newenvironment{proof}{{\it Proof. }}{\nopagebreak\hspace*{0.5cm}\hfill$\Box$\vspace{0.5cm}}

\def\bepf{\begin{proof}}
	\def\epf{\end{proof}}

\usepackage{algorithm}
\usepackage{algpseudocode}
\usepackage{adjustbox}
\usepackage{amsfonts} 

\usepackage{url}

\def\bfi#1{\textbf{#1}} 
\def\sca#1{{\sc #1}} 

\def\eeq{\end{equation}}
\def\lbeq#1{\begin{equation} \label{#1}}
\def\bary{\begin{array}}
	\def\eary{\end{array}}
\def\gzit#1{{\rm (\ref{#1})}} 			

\def\wh{\hat}
\def\wt{\tilde}
\def\pt{$\bullet$ }

\def\st {\fct{s.t.}}

\def\eps{\varepsilon}

\def\Rz{\mathbb{R}}

\def\fct#1{\mathop{\rm #1}}	                

\def\rank{\fct{rank}}

\def\argminx{\fct{argmin}}

\usepackage{array,multirow,graphicx}
\usepackage[utf8]{inputenc}

\usepackage{graphicx}

\usepackage{colortbl}

\usepackage[table]{xcolor}

\usepackage{fancyhdr}

\usepackage{array,multirow,graphicx}

\usepackage{etex}

\usepackage{longtable}
\usepackage{lscape}
\usepackage{hhline}
\usepackage{caption}
\usepackage{blindtext}
\usepackage{lmodern}
\usepackage{textcomp}
\usepackage{chngcntr}
\usepackage{tikz}
\usepackage{csquotes}

\fancypagestyle{plain}{\fancyhf{}\fancyfoot[C]{\bfseries \thepage}}

\newcolumntype{?}{!{\vrule width 1pt}}

\newcounter{subsubsubsection}[subsubsection]
\def\subsubsubsectionmark#1{}

\def\subsubsubsection{\@startsection
	{subsubsubsection}{4}{\z@} {-3.25ex plus -1
		ex minus -.2ex}{1.5ex plus .2ex}{\normalsize\bf}}
\def\l@subsubsubsection{\@dottedtocline{4}{4.8em}
	{4.2em}}

\newcolumntype{R}[2]{%
	>{\adjustbox{angle=#1,lap=\width-(#2)}\bgroup}%
	l%
	<{\egroup}%
}
\newcommand*\rot{\multicolumn{1}{R{30}{1em}}}

\usetikzlibrary{matrix,shapes,arrows,positioning,chains}

\def\bfi#1{{\bf{#1}}}
\def\sca#1{{\sc{#1}}} 



\tikzstyle{decision} = [diamond, draw, fill=blue!20,
text width=4.5em, text badly centered, node distance=3cm, inner sep=0pt]
\tikzstyle{block} = [rectangle, draw, fill=blue!20,
text width=3em, text centered, rounded corners, minimum height=2em]
\tikzstyle{line} = [draw, -latex']
\tikzstyle{cloud} = [draw, ellipse,fill=red!20, node distance=3cm,
minimum height=2em]

\usepackage{pgfplots}
\pgfplotsset{compat=1.18}
\usetikzlibrary{arrows.meta,plotmarks}

\begin{document}
	
	\begin{center}
		
		{\Large \bf A feasible point method for derivative-free }\\
		{\Large \bf linearly constrained optimization}\\

		\vspace{0.5cm}
		{\large \bf Morteza Kimiaei}
		\centerline{\sl Fakult\"at f\"ur Mathematik, Universit\"at Wien}
		\centerline{\sl Oskar-Morgenstern-Platz 1, A-1090 Wien, Austria}
		\centerline{\sl email: morteza.kimiaei@univie.ac.at}
		\centerline{\sl WWW: \url{http://www.mat.univie.ac.at/~kimiaei}}

		\vspace{0.5cm}
		{\large \bf Arnold Neumaier}
		\centerline{\sl Fakult\"at f\"ur Mathematik, Universit\"at Wien}
		\centerline{\sl Oskar-Morgenstern-Platz 1, A-1090 Wien, Austria}
		\centerline{\sl email: Arnold.Neumaier@univie.ac.at}
		\centerline{\sl WWW: \url{http://www.mat.univie.ac.at/~neum}}

	\end{center}

	{\bf Abstract.}
This paper presents a framework for linearly constrained derivative-free optimization (DFO) in which the original objective function is evaluated only at feasible points. The method first analyzes the linear constraint system using certificates of infeasibility and implied activity to detect infeasibility or a singleton feasible set, or, when the feasible set has positive dimension, to identify the affine subspace spanned by the feasible set and obtain an equivalent reduced representation with a strictly feasible base point. It then transforms the resulting reduced problem into an unconstrained optimization problem by means of a pointed feasible retraction and an exact penalty function. The retraction maps arbitrary trial points to feasible points by an explicit scaling in the reduced coordinates, so that no constrained optimization or projection subproblem is required during objective evaluations. This construction permits a broad range of unconstrained DFO methods to be used as inner solvers while preserving feasibility of the original objective evaluations. Numerical experiments compare several such variants on linearly constrained test problems of different dimensions. The results show that the choice of inner solver affects robustness, the number of problems solved, function-evaluation efficiency, and computational time, and that relative performance changes with problem dimension.

	{\bf Keywords.} Derivative-free optimization; linearly constrained optimization; penalty methods; benchmarking
	
	\vspace{0.2cm} {\em 2020 Mathematics Subject Classification: 90C56}.

	\hfill \today

	\vspace{1cm}
	\begin{sloppypar}
		
		\section{Introduction}
		In this paper, we consider the minimization of linearly constrained derivative-free optimization (DFO) problems. Throughout the paper, we denote by $x_I$ the subvector of a vector $x$ with indices taken from a list $I$ and by $A_{I:}$ the submatrix of a matrix $A$ with row indices taken from a list $I$. $|S|$ denotes the size of a set $S$. Thus, our optimization problems have 
		\lbeq{e.lcopt}
		\bary{ll}
		\mbox{min} & f(x)\\
		\st         & x\in C,
		\eary
		\eeq
		with the {\it feasible set}
		\lbeq{e.C}
		C:=\{x\in \Rz^n \mid (Ax)_E = b_E,~~~ (A x)_P\le b_P\}
		\eeq
		and the objective function $f$ is available only through computed function values at feasible points $x\in C$. Here, $A\in \Rz^{m\times n}$, $b\in \Rz^m$, and $(E,P)$ is a partition of $\{1,\ldots,m\}$. The theoretical analysis is stated for exact function values. In the numerical experiments, no artificial noise is added; the computed function values are nevertheless affected by the usual finite-precision and rounding errors. No prior information about gradients, Lipschitz constants, or the structure of $f$ is assumed to be available.
		
If $C$ is empty, the problem is infeasible. Otherwise, $C$ contains some feasible points $x^0$. In this case, we assume that $f:C\to\Rz$ is continuously differentiable on $C$, as a function on the affine subspace spanned by $C$. Derivatives are used
only in the analysis, and objective evaluations are required
only at feasible points. A {\it stationary point} of
\gzit{e.lcopt} is a point $\wh x\in C$ satisfying
$\partial_{x-\wh x} f(\wh x)\ge0$ for all $x\in C$, where $\partial_v$ denotes the directional derivative in the direction $v$. Every local minimizer
of \gzit{e.lcopt} is a stationary point.

This paper develops a solver, called {\it Vienna Linearly
	Constrained Derivative-Free Optimization} (\texttt{VLDFO}),
for solving \gzit{e.lcopt}. \texttt{VLDFO} is a
{\it feasible point method}: every evaluation of the original
objective function is performed at a feasible point. This property
is particularly important in applications where the objective
function is defined, reliable, or meaningful only on the feasible
set, so that evaluations at infeasible arguments are unavailable
or undesirable.

		{\bf Bounds on variables and two-sided linear constraints.}
		In the implementation, simple bounds on the variables are allowed, so that we may impose \(l\le x\le u\), where \(l\in(\Rz\cup\{-\infty\})^n\) and
		\(u\in(\Rz\cup\{+\infty\})^n\). We also allow two-sided linear constraints of the form $b^l_P\le (Ax)_P\le b^u_P$, where \(b^l_P\) and \(b^u_P\) denote the corresponding lower and upper right-hand sides. These constraints are not treated separately. Instead, all finite variable bounds and both sides of the linear constraints are rewritten as one-sided linear inequalities and incorporated into the inequality system.
		
		More precisely, let $I_l:=\{i\in\{1,\ldots,n\}\mid l_i>-\infty\}$ and $I_u:=\{i\in\{1,\ldots,n\}\mid u_i<+\infty\}$. Then, the finite lower and upper variable bounds are equivalently written as
		\(-x_{I_l}\le-l_{I_l}\) and \(x_{I_u}\le u_{I_u}\), respectively. Using the notation introduced above, if \(I\) denotes the \(n\times n\) identity matrix, these constraints become \(-I_{I_l:}x\le-l_{I_l}\) and \(I_{I_u:}x\le u_{I_u}\). Combining them with the two-sided linear constraints therefore gives the single inequality system
		\[
		\bar A_Px\le\bar b_P \ \ \mbox{with} \ \
		\bar A_P :=
		\pmatrix{
			A_{P:}\\
			-A_{P:}\\
			I_{I_u:}\\
			-I_{I_l:}
		}
		\ \ \mbox{and}\ \
		\bar b_P :=
		\pmatrix{
			b^{u}_{P}\\
			-b^{l}_{P}\\
			u_{I_u}\\
			-l_{I_l}
		}.
		\]
		If the original problem is already given in the one-sided form
		\((Ax)_P\le b_P\), then the block \(-A_{P:}x\le-b^l_P\) is simply omitted, or equivalently \(b^l_P=-\infty\). Thus, all supported variable bounds and two-sided linear constraints can be represented by the same one-sided inequality system. For simplicity of notation, throughout the rest of the paper we again write
		\(A_{P:}x\le b_P\) instead of \(\bar A_Px\le\bar b_P\), with the understanding that all finite variable bounds and applicable lower linear bounds have already been incorporated into the inequality index set \(P\).

		\bfi{Penalty Problem.} Define $a^+:=\max(0,a)$. If $f$ is defined or extended to $\Rz^n$, associated with the problem \gzit{e.lcopt} is the classical unconstrained penalty problem
		\lbeq{e.unopt0}
		\bary{ll}
		\mbox{min} & f_\gamma(x)=f(x)+\gamma\phi(\|r(x)\|)\\
		\st         & x\in\Rz^n,
		\eary
		\eeq
		with some norm $\|.\|$ and
		\[
		r_i(x):=\cases{(Ax-b)_i & if $i\in E$,\cr
			(Ax-b)_i^+ & if $i\in P$.}
		\]
		Here, $f_\gamma(x)$ is a {\it penalty function} with a penalty parameter $\gamma>0$, and the {\it penalizing functional} $\phi:\Rz_+\to\Rz_+$ is monotonically increasing and continuously differentiable with $\phi(0)=0$. Boundedness below of the penalty function depends on both the objective function and the penalty term and generally requires additional assumptions. If the constraints are violated, the penalty term $\gamma\phi(\|r(x)\|)$ is positive whenever $\phi(t)>0$ for $t>0$; otherwise it is zero. If the global minimizers of the penalty problem coincide with those of the original problem \gzit{e.lcopt}, the penalty function is called {\it exact}; otherwise {\it inexact}. When using the Euclidean norm and $\phi'(0)=0$, the resulting penalty function can be continuously differentiable under suitable regularity assumptions. Such smooth penalties are generally not exact for constrained problems with nonzero optimal constraint multipliers, although exactness may hold in particular cases.

\subsection{Overview}

\texttt{VLDFO} has two main phases:

\pt Phase~I analyzes the linear constraint system defining the feasible set \gzit{e.C} using certificates of infeasibility (COI) and
certificates for implied activities (CIA). A COI certifies infeasibility, whereas a CIA
identifies inequalities that are active at every feasible point;
see Proposition~\ref{p.implied} and Theorem~\ref{t.implied}.
A Phase~I procedure based on Theorem~\ref{t.main1}
successively identifies such implied activities and
incorporates them into the equality system. In this way, the
procedure identifies the affine subspace spanned by the feasible set $C$ and, if the latter has positive dimension, an interior point in this subspace. In the latter case, Theorem~\ref{t.main} then yields an equivalent reduced problem
for which the remaining inequalities admit the strictly feasible
base point \(z=0\). If $C$ does not have positive dimension, the problem is infeasible or has only a single feasible point, and the problem is solved at this stage.

\pt Phase~II uses this reduced representation
to transform the linearly constrained problem \gzit{e.lcopt} into an
unconstrained, generally nonsmooth penalty problem, defined by
\gzit{e.unprob} below. A pointed retraction maps each trial point
generated by the unconstrained inner solver to a feasible reduced point
by an explicit scaling (\gzit{e.betadef}, below), so that the original
objective is evaluated only at feasible points. Consequently, constraint
handling is separated from the optimization mechanism, allowing
local, global, deterministic, or randomized unconstrained DFO
methods to be used as inner solvers. Theorem~\ref{t.main2}
establishes the corresponding relations between the global and
local minimizers of the transformed and reduced constrained
problems.

In Phase~I, a single linearly constrained
problem suffices when no implied activity is detected, whereas
implied activities may require a finite sequence of such problems
until infeasibility is certified or the required reduced
representation is obtained. Thereafter, Phase~II requires no further
constrained optimization or projection subproblems, since every trial
point is mapped to feasibility directly by an explicit
scaling (\gzit{e.betadef}, below) in the reduced coordinates.

\subsection{Related work}

DFO algorithms are based on many years of experience with numerous
applications in science, engineering, industry, and chemistry; for
applications, see {\sc Audet \& Hare} \cite{AudH} and {\sc Conn}
et al. \cite{ConSV}; for classifications, see \sca{Larson} et al.
\cite{LarMW}; and for numerical behaviour, see \sca{Mor{\'e} \&
	Wild} \cite{MorW}, \sca{Rios \& Sahinidis} \cite{RioS},
\sca{Kimiaei} \cite{Kim}, and \sca{Kimiaei \& Neumaier}
\cite{KimN}. The principal DFO families considered here include
model-based methods \cite{ElsN0,GraTT,HuyN,Pow02,Pow08}, direct-search
methods \cite{Nomad,GraRVZ,Hig,HuyN0,PorT,VazV,VazVicente2009PSwarmLC},
line-search methods \cite{Kim,KimN,LucS}, and matrix-adaptation
evolution strategies \cite{AugH,Beyer2020,Beyer2017,MADFO}.

		{\bf Model-based methods} construct local models, often quadratic, by fitting or interpolation, solve these models subject to trust regions to generate trial points, and compare the actual reduction in the objective with the reduction predicted by the model. If the agreement is good, the trial point is accepted, and the trust region is maintained or expanded; otherwise, the trust region is reduced.
		
	{\bf Line search methods} considered here can be divided into two
	classes depending on whether approximate directional derivatives are
	used in the line-search acceptance conditions. In practice, methods
	that use approximate directional derivatives can perform poorly in
	the presence of noise because finite-difference estimates may become
	inaccurate. Line-search methods that do not use approximate
	directional derivatives instead base their acceptance conditions on
	function-value reductions and may employ extrapolation along
	successful search directions.
		
	{\bf Direct search methods} generate trial points along prescribed,
	adaptive, or randomized search directions and compare their objective
	values with that of the current point. An improving trial point may
	be accepted as the new iterate, with the step size maintained or
	increased, whereas unsuccessful polling typically leads to a
	reduction of the step size.
		
		{\bf Matrix adaptation evolution strategies} have three phases. In the first phase (the mutation phase), a finite number of distribution directions are generated, transformed by the adaptively determined affine scaling matrix into mutation directions, and the resulting mutation points are evaluated. In the second phase (the selection phase), the function values of the mutation points are sorted in ascending order. Then, the mutation directions and distribution directions are sorted accordingly. Then, some mutation points with the lowest function values and the corresponding distribution and mutation directions are selected, which are used in the next phase. In the third phase (the recombination phase), the recombination direction, which is a weighted average of the selected mutation directions, is calculated, and the recombination point, which is the sum of the old recombination point and the recombination direction scaled by a step size, is calculated. Both the affine scaling matrix and the step size are updated by using the selected information obtained from the mutation phase.
		
		\subsubsection{Solvers for unconstrained and bound constrained problems}\label{sec:un+bd}
		
	We first discuss DFO solvers that are used here in unconstrained or bound-constrained form to solve the nonsmooth penalty problems associated with \gzit{e.lcopt}; solver families considered through their linearly constrained variants are discussed instead in Subsection~\ref{sec:LC}.
		
		Model-based methods constitute a principal DFO class, with performance depending on interpolation geometry, dimension, noise, and model construction \cite{ConSV,Kim}. A full quadratic model has
		\[
		1+n+\frac{n(n+1)}2=\frac{(n+1)(n+2)}2
		\]
		coefficients, whereas underdetermined constructions require fewer evaluations. Powell's \texttt{UOBYQA} \cite{Pow02} uses full quadratic interpolation, while \texttt{NEWUOA} \cite{Pow08} typically uses $2n+1$ interpolation points and a minimum-change Hessian update; \texttt{BOBYQA} \cite{Powell2009BOBYQA} extends this philosophy to bound constraints. \texttt{BCDFO} \cite{GraTT} uses active-set trust-region models near bounds, while \texttt{DFOTR} \cite{Bandeira2012DFOTR} is another interpolation-based trust-region method. \texttt{SID-PSM} \cite{CusRV,CusV} combines pattern-search polling with quadratic models and simplex derivatives. Radial-basis models are used by \texttt{ORBIT} \cite{WildRegisShoemaker2008ORBIT}, while {\sc Curtis, Dezfulian \& W\"achter} \cite{CurtisDezfulianWachter2024} exploit surrogate information in structured bound-constrained DFO.

	Derivative-free line-search methods can be divided into
	two classes according to their acceptance conditions.
	The first class uses approximate directional derivatives,
	typically obtained by finite differences, in Armijo or
	Wolfe line-search conditions. Two representatives of this
	class are \texttt{fminunc} \cite{MatlabOTB} and
	\texttt{SSDFO} of {\sc Kimiaei, Neumaier \& Faramarzi}
	\cite{KimNF}. Although \texttt{fminunc} can therefore be called using function
	values only, it internally estimates derivatives and is used here as
	a comparison solver rather than as a genuine DFO method. The second class uses forcing functions based on function values instead of approximate directional
	derivatives in the acceptance conditions. This class
	includes the deterministic method \texttt{SDBOX} of
	{\sc Lucidi \& Sciandrone} \cite{LucS} and the randomized
	methods \texttt{VRBBO} of {\sc Kimiaei \& Neumaier}
	\cite{KimN} and \texttt{VRDFON} of {\sc Kimiaei} \cite{Kim}.
	These methods may employ extrapolation along successful
	search directions to accelerate convergence. The
	\texttt{DFLBOX} method of {\sc Liuzzi, Lucidi \& Rinaldi}
	\cite{Liuzzi2011} and the \texttt{DFNDFL} method of
	{\sc Giovannelli, Liuzzi, Lucidi \& Rinaldi}
	\cite{Giovannelli2022DFNDFL} also employ derivative-free
	line searches; the latter additionally uses nonsmooth
	continuous-search techniques and dense deterministic
	direction sequences.

		Classical direct-search methods include the Nelder--Mead simplex method \cite{NelderMead1965}, multidirectional search \cite{Torczon1989}, and related direct-search methods surveyed by \sca{Wright} \cite{wright1996direct}; in the present interfaces, \texttt{NELDER} is distinct from \texttt{NMSMAX}, while \texttt{ADSMAX} uses alternating coordinate searches \cite{Hig}. Bound-constrained global methods include \texttt{DIRECT} \cite{JonesPerttunenStuckman1993}, \texttt{MCS} \cite{HuyN0}, and the noise-oriented \texttt{SNOBFIT} \cite{HuyN}. \texttt{GLODS} \cite{CusM} combines global exploration with multiple local direct searches; in the current interface it operates within a finite auxiliary box.
		
		Evolution strategies include \texttt{CMAES}, whose covariance adaptation follows {\sc Hansen \& Ostermeier} \cite{HansenOstermeier2001} and whose implementation is cited through \cite{AugH}, together with \texttt{fMAES} \cite{Beyer2020}. {\sc Kimiaei \& Neumaier} proposed \texttt{MADFO} \cite{MADFO}, combining matrix adaptation with a stochastic nonmonotone line search and subspace mechanisms. The related \texttt{MATRS} \cite{MATRS} adds extrapolation and a model-based trust-region phase, while \texttt{DMADFO} \cite{DMADFO} replaces full affine scaling by diagonal scaling derived from conjugacy information. Together, these solvers provide complementary model-based, line-search, direct-search, global-search, and matrix-adaptation mechanisms used by \texttt{VLDFO} for the penalty problem associated with \gzit{e.lcopt}.

	\subsubsection{Solvers for linearly constrained problems}\label{sec:LC}

	 Dedicated DFO methods for linear constraints, however, employ several, often overlapping, mechanisms, including feasible direct search, model-based trust regions, randomized or evolutionary search, and projection- or mapping-based reformulations. Fletcher's equality-constrained penalty \cite{fletcher1970class} and related developments in exact penalties and nonlinear programming \cite{conntrust,di1984class,estrin2020implementing,huyer2003new,NocW} motivate such reformulations.
	
	A major direct-search lineage begins with {\sc Lewis \& Torczon} \cite{LewisTorczon2000LC}, who extended pattern search to linearly constrained minimization using feasible local search patterns adapted to the geometry of the feasible-set boundary. {\sc Kolda, Lewis \& Torczon} \cite{Kolda2007} developed corresponding generating-set stationarity results based on approximate tangent cones. Degenerate constraints and implementation issues were further treated by {\sc Lewis, Shepherd \& Torczon} \cite{LewisShepherdTorczon2007} and {\sc Abramson} et al. \cite{AbramsonBrezhnevaDennisPingel2008}, while active-set identification was studied by {\sc Lewis \& Torczon} \cite{LewisTorczon2009}. Related methods include \texttt{PSwarm} of {\sc Vaz \& Vicente} \cite{VazVicente2009PSwarmLC}, which combines particle-swarm exploration with pattern-search polling; \texttt{PSIFA} \cite{DinizEhrhardt2018,DinizEhrhardt2020}, which combines pattern-search ideas with implicit filtering for noisy linearly constrained problems; and the randomized feasible direct-search method \texttt{DSPFD} of {\sc Gratton} et al. \cite{GraRVZ19}.
	
	Model-based approaches include Powell's \texttt{COBYLA} \cite{COBYLA}, which uses linear interpolation models of the objective and constraints; \texttt{LCOBYQA} of {\sc Gumma} et al. \cite{Gumma2013}, which extends quadratic-interpolation ideas related to \texttt{NEWUOA} \cite{Pow08} to linearly constrained problems; and Powell's \texttt{LINCOA} \cite{LINCOA}, which combines quadratic models with an active-set treatment of linearly constrained trust-region subproblems. Other approaches include the interior-backtracking trust-region method of {\sc Gao \& Cao} \cite{GaoCao2018} and the probabilistic affine-scaling method of {\sc Zhao, Wang \& Zhu} \cite{ZhaoWangZhu2022}. Nonsmooth and specialized variants include the discrete-gradient method of {\sc Bagirov, Ghosh \& Webb} \cite{BagirovGhoshWebb2006}, the finite-minimax method of {\sc Liuzzi, Lucidi \& Sciandrone} \cite{LiuzziLucidiSciandrone2006}, and the nonmonotone direct-search method of {\sc Garc\'ia--Palomares} \cite{GarciaPalomares2020}.
	
	Global approaches include the linearly constrained \texttt{DIRECT}-type method of {\sc Stripinis \& Paulavi\v{c}ius} \cite{StripinisPaulavicius2023}. The globally convergent evolution strategies \texttt{GCES} of {\sc Diouane} et al. \cite{DioGV} use extreme barriers or projection in the constrained setting and may additionally exploit quadratic models. More general constrained solvers also cover the linear case, notably \texttt{MADS}/\texttt{NOMAD} \cite{Nomad,Audet2022NOMAD4,Dig11} and \texttt{SID-PSM} \cite{CusRV,CusV,SIDPSMmanual}. More recently, {\sc Royer, Sohab \& Vicente} \cite{RoyerSohabVicente2024} proposed a full--low framework combining projected finite-difference quasi-Newton steps with probabilistic direct search.
	
	To clarify their treatment of feasibility, we distinguish these linearly constrained DFO methods according to whether the original objective is evaluated only at feasible points and, when it is, how feasibility is enforced. This distinction concerns the points at which the original objective is actually evaluated, rather than merely the feasibility of accepted iterates or the algebraic generation of intermediate candidate points.
	
	(i) From the viewpoint of feasibility of the original objective evaluations, the feasible-point pattern-search lineage of {\sc Lewis \& Torczon} and {\sc Kolda, Lewis \& Torczon}, together with its subsequent extensions \cite{AbramsonBrezhnevaDennisPingel2008,Kolda2007,LewisShepherdTorczon2007,LewisTorczon2000LC,LewisTorczon2009}, adapts the search or polling directions to the local geometry of the feasible set and restricts admissible steps so that the objective is evaluated at feasible trial points. \texttt{PSwarm} \cite{VazVicente2009PSwarmLC}, \texttt{PSIFA} \cite{DinizEhrhardt2018,DinizEhrhardt2020}, and \texttt{DSPFD} \cite{GraRVZ19} likewise exploit the explicitly known linear constraints to keep objective evaluations within the feasible set. In particular, \texttt{PSwarm} may algebraically generate an infeasible polling candidate, but such a point is handled by an extreme barrier rather than by evaluating the original objective there, while its particle-swarm search step explicitly enforces feasibility before objective evaluation. The finite-minimax method of {\sc Liuzzi, Lucidi \& Sciandrone} \cite{LiuzziLucidiSciandrone2006} uses feasible search directions together with a maximum feasible step length, while the method of {\sc Garc\'ia--Palomares} \cite{GarciaPalomares2020} constructs a finite set of feasible search directions. In these direct-search approaches, feasibility of the objective evaluations is therefore enforced primarily through the known linear geometry and feasibility tests, rather than by solving a separate projection problem at every objective evaluation.
	
	(ii) Among model-based methods, \texttt{LCOBYQA} \cite{Gumma2013} starts from a strictly feasible point, constructs feasible interpolation points, and generates trust-region trial steps subject to the linear constraints; its geometry-improving procedures are likewise designed to preserve feasibility. This differs from \texttt{COBYLA} \cite{COBYLA}, which may evaluate infeasible points while balancing objective improvement and constraint violation, and from \texttt{LINCOA} \cite{LINCOA}, whose interpolation and geometry-management procedures may use points outside the feasible polyhedron. The interior-backtracking method of {\sc Gao \& Cao} \cite{GaoCao2018} and the affine-scaling approach of {\sc Zhao, Wang \& Zhu} \cite{ZhaoWangZhu2022} maintain feasible interior iterates through affine-scaling and backtracking mechanisms. However, feasibility of the iterates should not be confused with feasibility of all objective samples: their derivative-free models may be constructed from additional sampling points in a neighborhood of the current iterate that are not required to satisfy all linear inequalities. They are therefore not formulated as feasible-only objective interfaces in the same sense as the feasible direct-search methods above. Similarly, the discrete-gradient construction of {\sc Bagirov, Ghosh \& Webb} \cite{BagirovGhoshWebb2006} uses surrounding objective information and assumes that the objective is available on an ambient neighborhood; it is therefore not designed for an objective whose evaluation is defined only on the feasible polyhedron.
	
	(iii) The global and more general constrained methods exhibit similar differences. The linearly constrained \texttt{DIRECT}-type method of {\sc Stripinis \& Paulavi\v{c}ius} \cite{StripinisPaulavicius2023} employs a mapping intended to eliminate the infeasible region from the objective-calculation process, so that objective calculations are performed in the feasible region. The constrained \texttt{GCES} variants \cite{DioGV}
	either reject infeasible generated points through an extreme barrier
	or project generated points onto the feasible set before evaluating
	the objective. The constrained implementation of \texttt{SID-PSM} \cite{CusRV,CusV,SIDPSMmanual} likewise evaluates the objective only at feasible search and poll points; previously evaluated feasible points may then be reused to construct simplex derivatives or quadratic models. For the \texttt{MADS}/\texttt{NOMAD} family \cite{AudetDennis2004Filter,AudetDennis2009PB,Nomad,Audet2022NOMAD4,Dig11}, the situation depends on the selected constraint-handling mechanism. With explicitly known linear constraints, an extreme barrier can reject an infeasible trial point before the original objective is evaluated, whereas progressive-barrier and filter mechanisms deliberately allow infeasible intermediate points. Thus, feasible-only objective evaluations are not an intrinsic property of \texttt{MADS}/\texttt{NOMAD}, but depend on the constraint-handling strategy and on how the black-box evaluation interface is organized. The full--low method of {\sc Royer, Sohab \& Vicente} \cite{RoyerSohabVicente2024} constructs feasible projected or tangent-cone trial steps and maintains feasible accepted iterates, but its Full-Eval component forms finite-difference approximations from additional objective samples that are not, in general, guaranteed to satisfy all inequality constraints.
	
	(iv) Projection-based unconstrained reformulations provide another way to obtain feasible-only objective evaluations. {\sc Galvan, Sciandrone \& Lucidi} \cite{GalvanSciandroneLucidi2021} evaluate the original objective at the Euclidean projection of each trial point onto the feasible set; for a general polyhedral feasible set, computing this projection amounts to solving a convex quadratic projection problem. {\sc Norkin} \cite{Norkin2024} considers both Euclidean projection and a pointed projection defined using a prescribed interior feasible point, thereby providing an alternative to Euclidean projection while retaining objective evaluations in the feasible region. Thus, these reformulations preserve feasibility by applying a projection mechanism to the trial points generated during the unconstrained optimization process.

In contrast to the methods above, \texttt{VLDFO} first uses Phase~I to identify an equivalent reduced representation with a strictly feasible base point, requiring the solution of a single linear program when no implied activities are detected and only a small number of linear programs otherwise. Phase~II then applies an explicit reduced-coordinate retraction. Consequently, the original objective is evaluated only at validated feasible points, while the Phase~II solver sees the unconstrained nonsmooth penalty problem \gzit{e.unprob} associated with \gzit{e.lcopt}, without requiring the solution of further linear programs.

\section{Theory}\label{sec.Theory}
		
	In this section, we provide the relevant theory for transforming linearly constrained DFO problems into the unconstrained nonsmooth DFO problems, defined by \gzit{e.unprob} below. We first analyze the linear constraint system to determine whether it is infeasible or to identify constraints that are active at every feasible point. This analysis either certifies infeasibility, identifies a singleton feasible set, or, when the feasible set has positive dimension, determines the affine subspace spanned by the feasible set and a minimally active feasible point in this subspace. We then express the problem in reduced coordinates on this affine subspace so that the remaining inequalities admit a strictly feasible base point. Using this reduced representation, we introduce a retraction that maps arbitrary trial points to feasible points, ensuring that the original objective function is evaluated only at feasible arguments. Finally, we combine this retraction with a penalty term to obtain an unconstrained, generally nonsmooth problem and establish the relation between its minimizers and those of the original linearly constrained problem.
	
	We consider the linear constraint satisfaction problem (CSP)
		\lbeq{e.LCSPP}
		\bary{ll}
		\mbox{find} & x\in \Rz^n\\
		\st         & x\in C=\{x\in \Rz^n \mid (A x)_P\le b_P,\ \ (Ax)_E = b_E\}.
		\eary
		\eeq
		A feasible point of \gzit{e.LCSPP} is called {\it strictly feasible} 
		if it satisfies all inequality constraints with strict inequality,
		$(A x-b)_P<0$. This is a stronger condition than feasibility; there are 
		simple problems, where feasible points exist but none of them is strictly
		feasible. For example, if $u>0$ then the only feasible point of $x\ge 0$ and $u^Tx\le 0$
		is $x=0$, which satisfies all constraints with equality. The index $i\in P$ is called {\it active} if $(Ax)_i=b_i$ holds for a feasible point $x$. For a 
		consistent CSP \gzit{e.LCSPP}, an index $i\in P$ is called an 
		{\it implied activity} if it is active at every feasible point $x$. The index set 
		\lbeq{e.I}
		I:=\{i\in P:(Ax)_i=b_i\mbox{ for every }x\in C\}
		\eeq
		includes all implied activities, while the index set 
		\lbeq{e.Pc}
		P':=P\setminus I
		\eeq
		excludes all implied activities. A feasible point of \gzit{e.LCSPP} is called 
		{\it minimally active} if $(A x)_i<b_i$ for all $i\in P$ that are not 
		implied activities, i.e.,  either $P'$ is empty or $(Ax-b)_{P'}<0$.

		For $x,y\in \Rz^n$, we define $\inf(x,y)$ componentwise by
		$\inf(x,y)_i:=\min(x_i,y_i)$. We call $x,y\in\Rz^n$ {\it complementary} if any, and hence all, of the following equivalent conditions hold:\\
		(i) $\inf(x,y)=0$,\\
		(ii) $x\ge 0$, $y\ge 0$, $x_iy_i=0$ ($i=1,\cdots,n$),\\
		(iii) $x\ge 0$, $y\ge 0$, $x^Ty=0$. 
		We call $x,y$ {\it strictly complementary} if, in addition, there is 
		no index $i$ with $x_i=y_i=0$.

		\begin{prop}\label{p.implied}
			Suppose that $C$ is nonempty and $y\in\Rz^m$ satisfies
			\lbeq{e.certImp}
			A^Ty=0,~~~b^Ty=0,~~~y_P \ge 0.
			\eeq
			Then every index $i\in P$ with $y_i\ne 0$ is an implied activity.
		\end{prop}
		
		\bepf
		Any feasible point $x$ satisfies $y_i(b-Ax)_i=0$ for $i\in E$, hence
		\[
		\sum_{i\in P} y_i(b-Ax)_i
		=y^T(b-Ax)
		=b^Ty-(A^Ty)^Tx=0.
		\]
		Since each term in the sum is nonnegative, we find
		$\inf(y_P,(b-Ax)_P)=0$, implying that $(Ax)_i=b_i$ for all
		$i\in P$ with $y_i\ne 0$. Therefore any such $i$ gives an implied
		activity.
		\epf

		We call a vector $y$ with nonzero $y_P$ satisfying \gzit{e.certImp} a 
		{\it certificate for implied activities (CIA)}. If $y_i>0$ for all 
		implied activities $i\in P$, the certificate for implied activities is 
		called {\it maximal}.

		The feasible set $C=\{x\in \Rz^n \mid (Ax)_E = b_E,~ (A x)_P\le b_P\}$ is empty if and only if there is a vector $y$, called a {\it certificate of infeasibility} (COI) for $C$, satisfying 
		\lbeq{e.certInfeas}
		A^Ty=0, \ \ b^Ty<0, \ \ y_P\ge 0.
		\eeq

		The basis for our approach is the following result by 
		\sca{Goldman \& Tucker} \cite{Goldman19574TO}. 
		
		\begin{thm}\label{t.implied}
			Exactly one of the following holds:
			\\
			(i) A certificate of infeasibility exists.
			\\
			(ii) A strictly feasible point exists.
			\\
			(iii) A minimally active feasible point $x$ and a maximal certificate 
			$y$ for implied activities exist, and $y_P$, $(b-Ax)_P$ are strictly
			complementary.
		\end{thm}

		The knowledge of a minimally active feasible point allows one to write
		the feasible set in a particularly useful form.
		
		\begin{thm}\label{t.main}
			Let $x^0$ be a minimally active feasible point of \gzit{e.LCSPP} and put
			\lbeq{e.r0}
			E':=E\cup I,~~~B:=A_{{P'}\,:}S,~~~r^0:=(Ax^0-b)_{P'}
			\eeq
			where $S$ is a matrix whose $d$ columns form a basis for the null space of 
			$A_{E'\,:}$ and $I$ and $P'$ come from \gzit{e.I} and \gzit{e.Pc}, respectively. Then $r^0<0$ and the feasible set can be written as the 
			set of $x\in\Rz^n$ for which there is a vector $z\in\Rz^d$ such 
			that
			\lbeq{e.redLCSP}
			x=x^0-Sz,~~~Bz\ge r^0.
			\eeq
			This reduced linear CSP has the strictly feasible point $(x,z)=(x^0,0)$. 
		\end{thm}

		\bepf
		$x$ is feasible iff $(Ax-b)_E=0$, $(Ax-b)_I=0$, and $(Ax-b)_{P'}\le0$. Since $x^0$ is feasible, every $i\in I$ is an implied activity, hence
		the definition of $I$ and $P'$ gives 
		\[
		(Ax^0-b)_E=0,~~~(Ax^0-b)_I=0,~~~r^0=(Ax^0-b)_{P'}<0,
		\]
		with strict inequality. In particular, any feasible point $x$ satisfies
		\[
		A_{E'\,:}(x^0-x)=(Ax^0-Ax)_{E'}=(b-b)_{E'}=0.
		\]
		Therefore, $x^0-x$ is a linear combination of columns of $S$, by 
		definition of $S$. This shows that there is a vector $z$ with
		\lbeq{e.xSz}
		x=x^0-Sz.
		\eeq
		Conversely, the point \gzit{e.xSz} satisfies $(Ax)_{E'}=(Ax^0)_{E'}-A_{E'\,:}Sz=b_{E'}$.  Hence, it is feasible iff $(Ax-b)_{P'}\le 0$. Since $r^0-Bz=(Ax^0-b)_{P'}-A_{P'\,:}Sz=(Ax-b)_{P'}$, this is equivalent to $Bz\ge r^0$. 
		\epf
		
		If $x^0$ is strictly feasible, then it is minimally active and hence the index set $I$ is empty and Theorem \ref{t.main} applies.

		The following theorem gives the Phase~I alternatives used to obtain the vector $x^0$ required in Theorem \ref{t.main}.
		
		\begin{thm}\label{t.main1}
			Given a real value $\delta>0$, an arbitrary vector $x'$
			satisfying the equality constraints $(Ax')_E=b_E$, and
			a vector $u$ with
			\lbeq{e.defu}
			u_E=0, \ \ \delta u_P>((Ax')_P-b_P)_+,
			\eeq
			let $\pmatrix{x^0\cr\wh \eps}$ be a minimum point of
			the linear program (LP)
			\lbeq{e.phase1}
			\bary{lll}
			\wh \eps =&\min & \eps \\
			& \st & (Ax-u\eps)_E=b_E,~~~
			(Ax-u\eps)_P\le b_P,~~~ -\eps \le\delta.
			\eary
			\eeq
			At the solution, one of three cases happens:\\
			(i) $\wh\eps>0$. In this case, there is a certificate
			for infeasibility, a vector $y$ satisfying
			\gzit{e.certInfeas}.\\
			(ii) $\wh\eps=0$. In this case, $x^0$ has at least
			one implied activity.\\
			(iii) $\wh\eps<0$. In this case, $x^0$ is a strictly
			feasible point.
		\end{thm}
		
		\bepf
		Since $A\in\Rz^{m\times n}$, we define
		$P_+:=P\cup\{m+1\}$. By setting
		\[
		\wt A_{P_+}:=
		\pmatrix{A_{P:} & -u_P \cr 0 & -1},
		\quad
		\wt A_E:=\pmatrix{A_{E:} & 0},
		\quad
		x_{\eps}:=\pmatrix{x\cr\eps},
		\]
		\[
		b_{P_+}:=\pmatrix{b_P\cr\delta},
		\qquad
		c:=\pmatrix{0_n\cr1},
		\]
		we reformulate problem \gzit{e.phase1} as
		\lbeq{e.phase2}
		\bary{lll}
		\wh\eps=&\min&c^Tx_{\eps}\\
		&\st&
		\wt A_Ex_{\eps}=b_E,~~~
		\wt A_{P_+}x_{\eps}\le b_{P_+}.
		\eary
		\eeq
		
		Let $b_+:=\pmatrix{b_E\cr b_{P_+}}$.
		Because of \gzit{e.defu},
		$x'_{\delta}:=\pmatrix{x'\cr\delta}$
		is feasible for \gzit{e.phase2}.
		Since $\eps\ge-\delta$, the objective is
		bounded below. Hence the feasible linear
		program has a finite minimizer
		$x^0_{\wh\eps}:=\pmatrix{x^0\cr\wh\eps}$.
		Its feasibility gives
		\lbeq{e.eqfeas}
		\wt A_Ex^0_{\wh\eps}=b_E.
		\eeq
		
		The optimality conditions for \gzit{e.phase2}
		imply that there is an augmented multiplier vector
		\[
		y^+:=\pmatrix{y_E\cr y_{P_+}}
		=\pmatrix{y_E\cr y_P\cr y_{\wh\eps}},
		\qquad y:=\pmatrix{y_E\cr y_P}
		\]
		such that $\wt A^Ty^+=-c$,
		$y_{P_+}\ge0$,  and $\inf(y_{P_+},
		(b_+-\wt A x^0_{\wh\eps})_{P_+})=0$. The stationarity equation $\wt A^Ty^+=-c$
		gives
		\lbeq{e.eqzero}
		A_P^Ty_P+A_E^Ty_E=0,\qquad
		u_P^Ty_P+y_{\wh\eps}=1.
		\eeq
		
		The complementarity condition gives
		\lbeq{e.eqzero1}
		(y^+)^T(b_+-\wt A x^0_{\wh\eps})=0.
		\eeq
		
		Using \gzit{e.eqzero} and
		\gzit{e.eqzero1}, we obtain $(y^+)^Tb_+
		=(y^+)^T\wt A x^0_{\wh\eps}
		=-c^Tx^0_{\wh\eps}
		=-\wh\eps$. Consequently,
		\[
		b^Ty+\delta y_{\wh\eps}=-\wh\eps.
		\]
		If $\wh\eps\ge0$, the artificial constraint
		$-\eps\le\delta$ is strictly inactive
		because $\delta>0$. Complementarity
		therefore implies $y_{\wh\eps}=0$.  Hence, $A^Ty=0$, $b^Ty=-\wh\eps$, $y_P\ge0$,
		and $u_P^Ty_P=1$. If $\wh\eps>0$, then $b^Ty<0$.
		Thus $y$ satisfies \gzit{e.certInfeas}
		and is a certificate of infeasibility
		for the original CSP \gzit{e.LCSPP}.
		
		If $\wh\eps=0$, then $A^Ty=0$, $b^Ty=0$, $y_P\ge0$, and $u_P^Ty_P=1$. In particular, $y_P\ne0$. Hence, $y$ is a CIA satisfying \gzit{e.certImp}. By Proposition~\ref{p.implied}, $x^0$ has at least one implied activity.
		
		Finally, if $\wh\eps<0$, primal feasibility gives $(Ax^0-b)_P\le u_P\wh\eps<0$, since $u_P>0$. Moreover, $(Ax^0)_E=b_E$ because $u_E=0$. Thus $x^0$ is strictly feasible for the original CSP.
		\epf

		The representation \gzit{e.redLCSP} is useful since it reduces the linearly 
		constrained optimization problem \gzit{e.lcopt} to the equivalent problem
		\lbeq{e.redLCopt}
		\bary{ll}
		\min & e(z):=f(x^0-Sz)\\
		\st  & Bz\ge r^0
		\eary
		\eeq 
		with the strictly feasible point $z=0$. If $f$ is continuously differentiable with gradient $g$, the function $e$ is continuously differentiable with directional derivative $\partial_ue(z):=-\partial_{Su}f(x^0-Sz)$. Since $r^0<0$, the point $z=0$ is a strictly feasible starting point. 
		
		To ensure that $f$ is evaluated only at feasible points even when $z$
		is infeasible, we move the argument $z$ to $z':=z/\beta_z$, where
		\lbeq{e.betadef}
		\beta_z:=
		\max(\{1\}\cup\{ (Bz)_i/r^0_i \mid i\in P'\})\ge 1.
		\eeq
		Since $r^0<0$ and $\beta_z\ge (Bz)_i/r^0_i$ for all $i\in P'$, we have
		$Bz'\ge r^0$. Thus $x^0-Sz'$ is feasible, and $\beta_z=1$ iff
		$z$ is feasible. Inspired by penalty methods, we add a penalty term
		$\gamma \beta_z$ for some $\gamma>0$ to $e(z')$. This transforms the
		constrained DFO problem \gzit{e.redLCopt} into the unconstrained
		DFO problem
		\lbeq{e.unprob}
		\bary{ll}
		\min & e_{\gamma}(z):=e(z/\beta_z)+\gamma\beta_z.
		\eary
		\eeq
		
		Since $\beta_z$ is the maximum of finitely many affine functions
		and $\beta_z\ge1$, the mapping $z\mapsto z/\beta_z$ is locally
		Lipschitz continuous. Under our differentiability assumption on
		$f$, the transformed objective $e_\gamma$ is therefore locally
		Lipschitz continuous but generally nonsmooth. Unconstrained DFO
		methods suitable for nonsmooth objectives can thus be used as
		inner solvers for \gzit{e.unprob}. The exactness result below
		concerns minimizers and does not provide a convergence guarantee
		for a particular inner solver.
		
		The following result shows that any unconstrained minimizer of the
		penalty problem \gzit{e.unprob} is a minimizer of
		\gzit{e.redLCopt} and hence gives a minimizer of
		\gzit{e.lcopt}. Thus \gzit{e.unprob} is an exact penalty function.
		
		\begin{thm}\label{t.main2}
			Let $\mathcal F:=\{z\in\Rz^d \mid Bz\ge r^0\}$ and $\gamma>0$.\\
			(i) If $\wh z\in\mathcal F$ is a minimum point of problem
			\gzit{e.redLCopt}, then it is also a minimum point of problem
			\gzit{e.unprob}.\\
			(ii) If $z\notin\mathcal F$, then $\beta_z>1$,
			$z/\beta_z\in\mathcal F$, and
			$e_{\gamma}(z)>e_{\gamma}(z/\beta_z)$.\\
			(iii) If $\wh z\in\Rz^d$ is a minimum point of problem
			\gzit{e.unprob}, then $\wh z\in\mathcal F$ and it is also
			a minimum point of problem \gzit{e.redLCopt}.\\
			(iv) If $\wh z\in\mathcal F$ is a local minimizer of problem
			\gzit{e.unprob}, then it is a local minimizer of problem
			\gzit{e.redLCopt}. Conversely, if $\wh z$ is a local minimizer
			of problem \gzit{e.redLCopt} with $B\wh z>r^0$, then it is
			also a local minimizer of problem \gzit{e.unprob}.
		\end{thm}
		
		\bepf
		(i) Assume that $\wh z$ is a minimum point of \gzit{e.redLCopt}, i.e., $e(\wh z)\le e(z)$ for all $z\in\mathcal F$. Since $z/\beta_z\in\mathcal F$ for every $z\in\Rz^d$, $\beta_z\ge1$, and $\beta_{\wh z}=1$, we obtain
		\[
		e_{\gamma}(\wh z)=e(\wh z)+\gamma\le e(z/\beta_z)+\gamma\beta_z=e_{\gamma}(z)
		\]
		for all $z\in\Rz^d$. Thus $\wh z$ is a minimum point of \gzit{e.unprob}.\\
		(ii) If $z\notin\mathcal F$, then $(Bz)_i<r^0_i<0$ for some $i\in P'$. Hence $(Bz)_i/r^0_i>1$, so $\beta_z>1$. By \gzit{e.betadef}, $z':=z/\beta_z\in\mathcal F$, and therefore $\beta_{z'}=1$. Consequently,
		\[
		e_{\gamma}(z)=e(z')+\gamma\beta_z>e(z')+\gamma=e_{\gamma}(z').
		\]
		(iii) Assume that $\wh z$ is a minimum point of \gzit{e.unprob}. By (ii), $\wh z$ cannot be infeasible, since otherwise $e_{\gamma}(\wh z)>e_{\gamma}(\wh z/\beta_{\wh z})$, contradicting its minimality. Hence $\wh z\in\mathcal F$ and $\beta_{\wh z}=1$. For every $z\in\mathcal F$, we also have $\beta_z=1$, and therefore
		\[
		e(\wh z)+\gamma=e_{\gamma}(\wh z)\le e_{\gamma}(z)=e(z)+\gamma.
		\]
		Thus, $e(\wh z)\le e(z)$ for all $z\in\mathcal F$, proving that $\wh z$ is a minimum point of \gzit{e.redLCopt}.\\
		(iv) Assume that $\wh z\in\mathcal F$ is a local minimizer of \gzit{e.unprob}. Since $\beta_z=1$ for all $z\in\mathcal F$, we have $e_{\gamma}(z)=e(z)+\gamma$ on $\mathcal F$. Hence the local minimality of $e_{\gamma}$ at $\wh z$ implies the local minimality of $e$ over $\mathcal F$. Conversely, suppose that $\wh z$ is a local minimizer of \gzit{e.redLCopt} and $B\wh z>r^0$. Then a neighborhood of $\wh z$ is contained in $\mathcal F$, where $\beta_z=1$ and $e_{\gamma}(z)=e(z)+\gamma$. Therefore $\wh z$ is also a local minimizer of \gzit{e.unprob}.
		\epf
		
	Theorem~\ref{t.main2} concerns the exact objective $e$. In the
	numerical implementation, the value supplied by the
	objective oracle is the computed function value at the corresponding
	feasible argument. Thus, the exactness result established in
	Theorem~\ref{t.main2} concerns the transformed objective itself,
	whereas the convergence behavior depends on the properties of the
	selected unconstrained DFO inner solver. The theorem therefore establishes exactness
	of the deterministic transformation, but does not by itself give
	a convergence result for the selected unconstrained DFO method.
	Such convergence depends on the assumptions and mechanisms of the
	corresponding inner solver.

		\section{The \texttt{VLDFO} algorithm}
		\label{sec.vldfo.algorithm}
		
		The algorithm has two phases. Phase~I uses only the constraint data
		to identify a suitable feasible starting point and the equalities
		defining its feasible face. Phase~II applies an unconstrained DFO
		solver to \gzit{e.unprob}. The mathematical description uses the
		multiplier convention of \gzit{e.certImp} and \gzit{e.certInfeas}.
		The implementation uses the sign-reversed multiplier convention
		$\widetilde y=-y$ for its internal inequality multipliers. Hence, the theoretical COI conditions
		\(y_P\ge0\) and \(b^Ty<0\) are represented internally as
		\(\widetilde y_P\le0\) and \(b^T\widetilde y>0\); the two conventions
		are mathematically equivalent.
		
		\medskip
		\noindent\textbf{High-level algorithm (\texttt{VLDFO}).}
		
		\noindent\textbf{Input:} The data of \gzit{e.lcopt}, a reference
		point, a Phase~I parameter $\delta>0$, a penalty parameter
		$\gamma>0$, an inner unconstrained DFO solver, and numerical stopping tolerances
		and budgets.
		
		\noindent\textbf{Output:} A feasible evaluated point with
		termination information, a validated COI, or an unresolved
		preprocessing outcome.
		
		\begingroup
		\renewcommand{\theenumi}{S\arabic{enumi}}
		\renewcommand{\labelenumi}{\theenumi.}
		\begin{enumerate}
			\item \label{it.vldfo.preprocess}
			Incorporate the finite constraint sides and variable bounds into
			the representation used in \gzit{e.LCSPP}. Check the equality
			subsystem and obtain an equality-feasible reference point.
			A multiplier obtained from this equality-consistency check is retained
			as a diagnostic certificate candidate; terminal infeasibility is
			decided subsequently through the Phase~I certificate-validation route.
			
			\item\label{it.vldfo.phase1}
			For the current partition of the rows, choose the relaxation
			vector according to \gzit{e.defu}. Obtain and validate primal
			candidates and multipliers for the LP problem \gzit{e.phase1}, distinguishing
			regularized primal recovery from LP optimality verification.
			
			\item\label{it.vldfo.faces}
			Apply the alternatives of Theorem~\ref{t.main1}. A validated COI
			terminates preprocessing. A strictly feasible witness completes
			Phase~I. In the CIA case, promote the certified implied-active
			inequalities to equalities and return to
			\ref{it.vldfo.phase1}. If the required numerical checks remain
			inconclusive, report an unresolved outcome. A validated singleton
			face can terminate the optimization directly.
			
			\item\label{it.vldfo.reduction}
			Construct the null-space basis and reduced data specified in
			\gzit{e.r0}. If the reduced dimension is zero, evaluate the unique
			feasible point and terminate. Otherwise, initialize the reduced
			iteration at the origin.
			
			\item\label{it.vldfo.oracle}
			For every inner trial vector, compute the scaling $\beta_z$ from
			\gzit{e.betadef}, retract the reduced argument, and recover the
			original argument through \gzit{e.xSz}. Validate this argument
			before calling the original objective oracle, and supply the
			inner solver with the value prescribed by the unconstrained nonsmooth problem \gzit{e.unprob}.
			
			\item\label{it.vldfo.inner}
			Run the selected unconstrained DFO method using this evaluation
			interface until a stopping condition is reached. Record original
			objective evaluations separately from transformed objective
			values.
			
			\item\label{it.vldfo.output}
			Apply the same retraction and validation to the final inner
			output. Return its feasible original-space counterpart and its
			original objective value, retaining the best evaluated feasible
			point and the evaluation and timing records separately.
		\end{enumerate}
		\endgroup
		
		\subsection{Constraint preprocessing and equality consistency}
		\label{subsec.vldfo.preprocess}
		
		\ref{it.vldfo.preprocess} uses the treatment of variable bounds
		and two-sided rows described in the introduction. Only finite
		constraint sides are retained. Positive row scaling may improve
		conditioning without changing the feasible set in \gzit{e.C}.
		The original unscaled data are retained for independent
		feasibility checks.
		
		For benchmark inputs, the default convention treats upper sides
		at or above $10^{20}$ and lower sides at or below $-10^{20}$ as
		representations of infinity. This convention can be disabled
		when all finite sides are to be retained.
		
		The equality subsystem must be checked before using
		Theorem~\ref{t.main1}. To describe the consistency test, let
		$D_E$ be a positive diagonal row-scaling matrix, and write
		$\widehat A_E=D_E^{-1}A_{E:}$ and
		$\widehat b_E=D_E^{-1}b_E$. From a reference point
		$x^{\rm ref}$, compute
		\[
		x^{\rm eq}
		=x^{\rm ref}+\widehat A_E^{\dagger}
		(\widehat b_E-\widehat A_E x^{\rm ref}),
		\qquad
		r_E=\widehat b_E-\widehat A_E x^{\rm eq},
		\]
		where $\dagger$ denotes the Moore--Penrose inverse. In exact
		arithmetic, a vanishing residual provides the reference point
		required by Theorem~\ref{t.main1}. If the residual is nonzero,
		the corresponding residual multiplier provides a natural
		Farkas-type candidate for inconsistency of the equality subsystem.
		In floating-point arithmetic this multiplier is checked using the
		same scale-aware certificate conditions employed in Phase~I.
		The implementation retains such a multiplier as diagnostic
		information; terminal infeasibility is not declared from this
		preliminary equality projection alone, but is decided by the
		subsequent Phase~I certificate-validation route. If there are no
		equality rows, the supplied reference point can be used directly.
		
		\subsection{Phase~I solves and numerical validation}
		\label{subsec.vldfo.phase1}
		
		\ref{it.vldfo.phase1} is applied to the current row partition.
		For each remaining inequality, the implementation chooses
		\lbeq{e.vldfo.uinit}
		u_i=1+\frac{2}{\delta}
		(a_i^T x^{\rm eq}-b_i)_+,
		\eeq
		where $a_i^T$ is the corresponding row of the current matrix.
		The equality components are set as prescribed in
		\gzit{e.defu}. This construction satisfies \gzit{e.defu}, including
		when $\delta$ is small, and provides the feasible augmented
		starting pair used in the proof of Theorem~\ref{t.main1}.
		After a successful face update, the retained feasible witness
		serves as the next reference point.
		
		The default initial value is $\delta_0=10^{-3}$, and the
		default primal feasibility tolerance is $10^{-5}$. Certificate
		verification uses separate, stricter tests.
		
		The three Phase~I solver routes are \texttt{minq8sep}~\cite{HuyerNeumaier2018MINQ8},
		\texttt{fmincon}~\cite{MathWorksFmincon}, and
		\texttt{quadprog}~\cite{MathWorksQuadprog}. The first two are used
		primarily for regularized primal recovery. With the notation of
		\gzit{e.phase2}, their objective is
		\lbeq{e.vldfo.qrho}
		q_\rho(x_{\varepsilon})
		=c^T x_{\varepsilon}
		+\frac{\rho}{2}\|x_{\varepsilon}\|_2^2,
		\qquad \rho=10^{-8},
		\eeq
		under the constraints of \gzit{e.phase2}. A Newton projection
		can also recover a primal candidate when the primary
		\texttt{minq8sep} return is not directly usable.
		
		In the implementation used for the numerical experiments, \texttt{minq8sep} solves the linearly constrained separable convex quadratic problem by an active-set-type method, while \texttt{fmincon} is used with its active-set algorithm. The primary \texttt{quadprog} call uses the interior-point-convex algorithm; if this call fails numerically, the implementation retries it and may finally use the active-set algorithm as a fallback.
		
		The numerical routes are invoked hierarchically. A validated
		strictly feasible \texttt{minq8sep} return terminates the
		fallback chain. Otherwise the conditional \texttt{fmincon}
		challenge is invoked. A validated strictly feasible
		\texttt{fmincon} return is again decisive. If its Phase~I value
		is numerically zero, the method proceeds directly to the CIA
		and facial-reduction branch only when a usable CIA route is
		already available. Otherwise the final \texttt{quadprog} route
		uses a zero Hessian to solve the unregularized Phase~I LP
		\gzit{e.phase1}. Thus the regularized solvers primarily provide
		robust primal recovery, whereas the unregularized
		\texttt{quadprog} route supplies the authoritative LP fallback
		when the preceding information is insufficient.
		
		Every primal candidate is checked independently against
		\gzit{e.phase1}. Multiplier reliability additionally requires
		the appropriate stationarity, complementarity, multiplier-sign,
		and primal-feasibility tests. When relating native solver
		multipliers to the analysis, their sign and row-scaling
		conventions must be accounted for when applying
		\gzit{e.certImp} or \gzit{e.certInfeas}. Newton-estimated
		multipliers are diagnostic candidates rather than accepted
		certificates. In particular, a positive regularized-QP
		objective or a positive candidate value of $\varepsilon$ does
		not by itself establish infeasibility. A directly validated
		strictly feasible point, however, is a sufficient feasibility
		witness without an optimality claim.
		
		A solver exit flag or a small unscaled residual alone is not
		used as a certificate. For a candidate COI multiplier $y$,
		\texttt{VLDFO} first applies a scale-aware direct test of the
		conditions
		\lbeq{e.vldfo.coi}
		A^Ty\approx0,\qquad b^Ty>0,
		\qquad y_P\leq0,
		\eeq
		together with nonzero and dimension checks. Because a useful
		Farkas multiplier may nevertheless be poorly scaled,
		\texttt{VLDFO} also evaluates the normalized diagnostic
		\lbeq{e.getI}
		f_I:={f_I^N}/{f_I^D},
		\qquad
		f_I^N:=\|A^Ty\|_2\,b_{\max},
		\qquad
		f_I^D:=A_{\max} \,(b^Ty).
		\eeq
		Here,
		\lbeq{e.getIscales}
		A_{\max}:=\max\{1,\max_{i,j}|A_{ij}|\},
		\qquad
		b_{\max}:=\max\{1,\max_i|b_i|\},
		\eeq
		with excessively large artificial right-hand sides excluded from
		the calculation of $\max b$. Before evaluating \texttt{getI}, the multiplier is positively
		normalized and its inequality components are checked for the
		required sign. The zero multiplier is rejected, since
		$A^Ty=0$ and $b^Ty=0$ for $y=0$ carry no information about
		infeasibility. A finite value satisfying
		\lbeq{e.getItest}
		f_I<1,\qquad f_I^D>0
		\eeq
		passes this normalized certificate test. The quantity $f_I$
		therefore measures the stationarity residual relative to the
		positive Farkas separation $b^Ty$, while accounting for the
		scaling of $A$ and $b$.
		
		Intermediate \texttt{getI} results obtained from regularized
		solver multipliers are retained as diagnostics or certificate
		candidates and do not by themselves replace the unregularized
		Phase~I LP route. For the final multiplier obtained through the
		unregularized Phase~I LP route, \texttt{VLDFO} evaluates both the
		direct scale-aware COI validation associated with
		\gzit{e.vldfo.coi} and the normalized \texttt{getI} test
		\gzit{e.getI}--\gzit{e.getItest}. The latter provides an additional
		normalized diagnostic for the selected certificate. If the required
		COI validation cannot be established, the infeasibility
		classification remains unverified. The multiplier, constraint
		system, and validation data are always kept consistent with the
		same current row partition.

		Numerical sign classification uses a tolerance fixed from the
		initial parameter $\delta_0>0$:
\[		
		\tau_{\varepsilon}
		=\min\bigl\{0.1\delta_0,
		10\sqrt{\epsilon_{\rm mach}}(1+\delta_0)\bigr\}.
\]
		Primal feasibility and certificate identities have separate
		tolerances. Values within the sign tolerance initiate the CIA
		verification branch; they do not establish an exact zero
		optimum by themselves. If a numerical retry enlarges $\delta$,
		the implementation uses $\delta\leftarrow\min\{1.1\delta,10^{10}\}$, 
		while retaining the original sign tolerance. A failed
		verification remains unresolved and is not reported as a
		certificate.
		
		The reported outcomes preserve this distinction:
		\texttt{COI} certifies infeasibility;
		\texttt{CIA unverified} and \texttt{COI unverified} mean that
		the required certificate could not be validated; and
		\texttt{Phase~I failed} means that preprocessing did not produce
		a validated state suitable for Phase~II. An unresolved or failed
		preprocessing run does not authorize an objective evaluation at
		an unvalidated point.
		
		In the experiments reported in Section~\ref{s.num}, no \texttt{CIA unverified}, \texttt{COI unverified}, or \texttt{Phase~I failed} outcome occurred. This empirical observation does not constitute a guarantee in finite-precision arithmetic; such outcomes may occur for other test collections, particularly when the required numerical validations cannot be established reliably.
		
		\subsection{Repeated facial reduction}
		\label{subsec.vldfo.faces}
		
		\ref{it.vldfo.faces} supplies the repeated preprocessing
		needed to obtain the minimally active point required by
		Theorem~\ref{t.main}. Let $E_j$ and $P_j$ denote the equality and
		inequality row sets at preprocessing round $j$, initially the
		sets in \gzit{e.LCSPP}. All certificate checks at round $j$ use
		this current partition.
		
		In case~(i) of Theorem~\ref{t.main1}, preprocessing terminates
		only after a COI has passed \gzit{e.certInfeas}. In case~(iii),
		the original feasibility and strict remaining slack of the
		primal witness are checked before accepting it as $x^0$.
		Case~(ii) requires another reduction round. For a validated CIA
		$y^{[j]}$, take
		\[
		J_j=\{i\in P_j:y_i^{[j]}>0\},
		\qquad
		E_{j+1}=E_j\cup J_j,
		\qquad
		P_{j+1}=P_j\setminus J_j.
		\]
		The normalization following \gzit{e.eqzero} in the proof of
		Theorem~\ref{t.main1} ensures that this support is nonempty in
		the exact zero-optimum case. Proposition~\ref{p.implied}
		justifies every promotion, so the feasible set is unchanged.
		A row is not promoted merely because it is active at the
		current primal point.
		
		The implementation can combine several independently validated
		CIA candidates on the same face. Positive combinations retain
		their certificate properties, and their inequality supports
		cannot cancel. Numerical support cutoffs are applied only after
		the multiplier itself has passed the corresponding CIA checks.
		If the initial support cutoff produces no usable row, the cutoff
		may be reduced in a scale-aware manner, and an auxiliary CIA
		search may provide additional candidate implied active rows.
		
		Every proposed move from $P_j$ to $E_j$ subsequently passes an
		independent implication-verification gate. For a candidate set
		$J\subseteq P_j$, the implementation considers the maximum
		attainable sum of its feasible slacks,
		\[
		\sup_{x\in C}
		\sum_{i\in J}(b_i-a_i^Tx).
		\]
		Since every feasible slack is nonnegative, a validated optimum
		equal to zero establishes that every row in $J$ is implied
		active. Candidate rows are first verified in batches; when a
		batch cannot be certified, it is split recursively until
		subsets or individual rows are classified as verified,
		disproved by a feasible positive-slack witness, or numerically
		unresolved. Only verified rows are promoted. Hence a failed
		LP or certificate calculation does not constitute mathematical
		evidence against implied activity.
		
		This mandatory implication-verification gate is distinct from
		the optional maximum-slack fallback used when the ordinary
		solver CIA, auxiliary CIA, and adaptive CIA mechanisms produce
		no candidate support. In that fallback, near-active rows are
		tested individually through their maximum attainable slack,
		\[
		\sigma_i=\sup_{x\in C}(b_i-a_i^Tx),
		\qquad i\in P_j.
		\]
		
		Before committing a verified promotion, the retained feasible
		witness is also checked to ensure that the proposed new
		equalities are numerically active. An auxiliary feasibility
		solve can recover such a witness when necessary. If no adequate
		witness can be obtained, the move is deferred rather than
		declared mathematically invalid. The implementation also retains
		the previously verified face before a promotion so that a
		numerically unsuccessful subsequent Phase~I solve can roll back
		the tentative face update.
		
		Each successful zero-optimum round removes at least one
		inequality. Hence there are at most as many such rounds as
		initial inequality rows in exact arithmetic. Numerical retries
		are distinct from these successful reductions. At a termination
		with strict remaining slack, all implied activities have been
		identified, and the final row sets agree with
		\gzit{e.I}, \gzit{e.Pc}, and \gzit{e.r0}.
		If the current equality matrix instead reaches full column rank,
		its unique point is tested against all original constraints.
		When this test succeeds, the feasible set is a singleton and
		the remaining implied rows need not be enumerated.
		
		If no inequalities remain, facial reduction is also complete:
		optimization continues in the subspace defined by the equalities,
		unless the feasible set is a singleton.
		
		\begin{remark}
			\label{rem.multipleCIA}
			The facial reduction step does not require all implied-active
			inequalities identified in one round to belong to the support of a
			single CIA. Suppose that, for the current partition $(E_j,P_j)$,
			$y^{(1)},\ldots,y^{(q)}$ are independently validated CIAs, so that
			\[
			A^T y^{(k)}=0,\qquad
			b^T y^{(k)}=0,\qquad
			y^{(k)}_{P_j}\ge0,
			\qquad k=1,\ldots,q.
			\]
			Then Proposition~\ref{p.implied} implies that every index
			$i\in P_j$ for which $y_i^{(k)}>0$ for at least one $k$ is an
			implied activity. Hence the union
			\[
			J_j
			:=
			\bigcup_{k=1}^q
			\{i\in P_j:y_i^{(k)}>0\}
			\]
			may be moved from $P_j$ to $E_j$ in the same facial-reduction round.
			Indeed, for any coefficients $\alpha_k>0$, the positive combination
			\[
			\bar y:=\sum_{k=1}^q \alpha_k y^{(k)}
			\]
			is again a CIA, and the nonnegativity of the inequality multipliers prevents cancellation of their positive supports. Thus it is not necessary to construct one certificate whose support contains all
			indices in $J_j$ simultaneously.
			
			This observation is also important for the numerical implementation.
			The regularized Phase~I solves used for primal recovery need not solve
			the exact linear program~\gzit{e.phase1}. Consequently, a computed
			value of $\eps$ lying within the numerical zero tolerance does not, by
			itself, prove that the exact Phase~I optimum satisfies
			$\widehat\eps=0$. Such a value is used only to initiate the CIA
			verification branch. A row is promoted to the equality system only
			after the current original constraint system has a validated feasible
			point and an independently validated CIA certifies that row as implied
			active. Different promoted rows may therefore be justified by
			different exact CIA certificates.
			
			Accordingly, once at least one nonempty set of implied-active rows has
			been certified, facial reduction may proceed immediately without
			requiring another Phase~I solver merely to enlarge that set. Additional
			CIAs may identify more implied activities in the same round and can
			improve computational efficiency, but they are not required for the
			validity of the reduction. Any remaining implied activities can be
			identified in subsequent facial-reduction rounds.
		\end{remark}
		
		\subsection{Construction of the reduced coordinates}
		\label{subsec.vldfo.reduction}
		
		\ref{it.vldfo.reduction} is performed after the preceding
		checks have supplied the point required by
		Theorem~\ref{t.main}. Construct the data through \gzit{e.r0},
		choosing an orthonormal null-space basis by a rank-revealing
		singular-value decomposition. Its dimension is $d=n-\rank(A_{E':})$. Positive row normalization and a consistent numerical rank
		tolerance are used in the face-rank test and basis construction.
		For a wide equality matrix, the full right singular-vector
		matrix is needed to retain every null-space direction.
		
		When $d=0$, the validated feasible point is unique, so one
		original objective evaluation completes the optimization.
		When $d>0$, initialize at $z_0=0$ and use \gzit{e.redLCopt}.
		Every remaining component of the residual in \gzit{e.r0} must
		have strict slack before the radial map is enabled.
		An empty inequality set does not eliminate the equality
		constraints: the null-space construction is still required.
		The identity basis is appropriate when the equality set is empty.
		
		\subsection{The feasible evaluation interface}
		\label{subsec.vldfo.oracle}
		
		In \ref{it.vldfo.oracle}, an unrestricted inner trial $z$
		first determines $\beta_z$ through \gzit{e.betadef}. Denote the
		retracted reduced argument by $\bar z=z/\beta_z$, and obtain
		the associated original argument from \gzit{e.xSz} with
		$z$ replaced by $\bar z$. Theorem~\ref{t.main} and the
		construction preceding Theorem~\ref{t.main2} establish its
		feasibility. The original oracle is called at this recovered
		argument; the penalty is added afterwards according to
		\gzit{e.unprob}.
		
		The recovered point may lie on the boundary, even though the
		base point has strict slack for the remaining inequalities.
		
		The real-valued scaling in \gzit{e.betadef} is used throughout;
		there is no integer rounding. Its empty-set convention also
		handles the absence of remaining inequalities. A finite positive
		$\gamma$ can be prescribed, with default value $1$ in the
		implementation. An optional scale initialization uses three
		independent standard Gaussian vectors $\xi_j\in\Rz^d$.
		Query the evaluation interface with unit penalty parameter
		at the origin and at $0.1\xi_j$. If $w_0,w_1,w_2,w_3$ are
		the returned values, set
		\[
		\gamma=\max\left\{10^{-2},
		\max_{1\le j\le3}|w_j-w_0|\right\}.
		\]
		All objective calls made by this initialization
		use the same feasible evaluation interface.
		
		Floating-point safeguards are applied separately from this
		mathematical definition. Nonfinite or complex trial data are
		rejected. The computed original argument is checked against
		the original unscaled constraints, including variable bounds,
		as well as the current reduced representation. If the reduced
		retraction fails its numerical feasibility test, the routine
		can contract the reduced step towards the origin, using factors
		$2^{-j}$ and at most $50$ contractions, and can fall back to
		the validated base point. A trial that fails the final checks
		receives a numerical rejection value without an original
		objective call. Numerical feasibility refers to the specified
		primal tolerance; it is distinct from exact-arithmetic membership
		in \gzit{e.C}.
		
		\subsection{Unconstrained inner optimization}
		\label{subsec.vldfo.inner}
		
		\ref{it.vldfo.inner} applies the chosen routine to
		\gzit{e.unprob} in unconstrained reduced coordinates. Every
		trial, including points outside the reduced feasible region,
		passes through \ref{it.vldfo.oracle}. Thus, constraint
		handling belongs to the evaluation interface, and the inner
		routine does not need to enforce the original linear constraints.
		The same interface is used for local, global, deterministic,
		and randomized inner methods.

		As explained before Theorem~\ref{t.main2}, the transformed
		objective is generally nonsmooth. It is locally Lipschitz when
		the original objective is locally Lipschitz on its feasible
		domain. A routine suitable for nonsmooth DFO may therefore be selected
		from the inner methods discussed with Table~\ref{a.sol} in
		Appendix~\ref{sec.inner-solvers}. Applying a routine designed for smooth
		objectives is an empirical algorithmic choice; its smooth convergence
		theory does not automatically cover the transformed problem. The
		distinction between numerical applicability and theoretical convergence
		for the generally nonsmooth transformed problem is discussed further in
		Appendix~\ref{sec.nonsmooth-applicability}. The evaluation interface
		itself requires objective values and linear algebra, without derivative
		or subgradient approximations.

		Theorem~\ref{t.main2} interprets minimizers of the transformed
		problem. It does not assert that a finite run of an arbitrary
		inner routine obtains such a minimizer. Termination by an
		evaluation budget, a time budget, or an inner stopping rule
		therefore carries its own recorded termination reason.

		\subsection{Returned point, stopping, and performance records}
		\label{subsec.vldfo.output}
		
		\ref{it.vldfo.output} applies the same retraction to the
		final reduced output and recovers its original counterpart
		through \gzit{e.xSz}. This point is validated before its
		original objective value is returned. If final recovery fails,
		the implementation uses the base point only after validating
		that point as well.
		
		The best feasible evaluation is a separate record. If
		$x^{(1)},\ldots,x^{(k)}$ are the valid original arguments
		evaluated so far and $v_1,\ldots,v_k$ are the values returned by the original oracle, define
		\[
		j_k\in\argminx_{1\le j\le k}v_j,
		\qquad
		x_k^{\rm best}=x^{(j_k)},
		\qquad
		v_k^{\rm best}=v_{j_k}.
		\]
		These records use the original oracle responses before adding
		the penalty. A finite inner run can return a point different
		from $x_k^{\rm best}$, so the returned point, the best feasible
		record, and the transformed objective value must remain
		distinguishable in the numerical results.
		
		The function-evaluation budget counts actual original oracle
		calls, including optional scale-initialization calls and a
		final evaluation when performed. Linear-algebra operations,
		Phase~I solver calls, and rejected trial vectors do not count
		as original objective evaluations. Phase~I solver statistics
		record calls actually made; a validated primal candidate is
		recorded separately from a verified LP optimum or certificate.
		
		A point obtained by Newton recovery is recorded as a recovered
		candidate and does not turn an unsuccessful \texttt{minq8sep}
		call into a successful one. Each solver's multiplier and
		validation diagnostics are retained together.
		
		For timing comparisons, the charged CPU time is
		$T_{\rm total}=T_{\rm I}+T_{\rm II}$. A completed feasible
		Phase~I state may be reused across different inner
		\texttt{DFO} solvers only after the stored state is checked
		against the current problem data and the current original
		constraints. The stored implied equalities and the strict-slack
		condition required by Phase~II are revalidated; feasibility of
		the cached point alone is insufficient. Reuse retains the
		stored Phase~I CPU cost in $T_{\rm total}$, whereas cache
		input/output overhead is excluded. Current-call solver counts
		are reported separately from the historical counts stored with
		the cached state.
		
		A terminal infeasibility result is stored separately from a
		feasible Phase~I state. Such an infeasibility cache contains
		the multiplier and the corresponding certificate system and is
		accepted on a later call only after the stored certificate has
		again passed the current validation procedure. In particular, the cached multiplier must again pass the current scale-aware COI validation; the normalized \texttt{getI} test is recomputed as an additional diagnostic.
		
		\section{Numerical results}\label{s.num}

		In this section, using the 246 linearly constrained \texttt{CUTEst} test problems~\cite{CUTEst}, whose characteristics are provided in the accompanying GitHub repository~\cite{VLDFOsupp}, we compare the 26 \texttt{VLDFO} variants with three linearly constrained DFO solvers (\texttt{DSPFD}, \texttt{LINCOA}, \texttt{NOMAD-4.6}) and one nonlinearly constrained DFO solver (\texttt{COBYLA}) summarized in Table~\ref{a.sol} of Appendix~\ref{sec.inner-solvers}.
		
		Since the present study focuses on linearly constrained problems, we emphasize dedicated linearly constrained solvers and include \texttt{COBYLA} only as a representative general nonlinear constrained DFO method; a broader comparison with nonlinear constrained solvers is beyond the scope of this study, whose main purpose is to assess the proposed theory and to investigate which unconstrained inner solvers provide the most robust and efficient \texttt{VLDFO} variants.

		The remainder of this section is organized as follows. We first define the solved-problem criteria, cost measures, efficiency scores, anomaly counters, and solver-selection strategy. We then examine the Phase~I behavior of \texttt{VLDFO}, including infeasibility detection and problems with singleton feasible sets. Next, we report the numerical comparisons for dimensions $n\in[1,50]$ and for the complete range $n\in[1,500]$, including the larger-dimensional subset $n\in[51,500]$, solved-case diagnostics, and the problems unsolved by all tested methods. Finally, we compare the two \texttt{VLDFO-NOMAD} variants with the direct linearly constrained \texttt{NOMAD-4.6} solver for dimensions $n\in[1,10]$.
		
		\subsection{Benchmarking framework and solver selection}
		\label{s.benchmark-framework}
		
		\bfi{Solved and unsolved problems.} We say that a solver \(s \in \mathcal{S}\) solves a problem (i.e., \bfi{solved}) if at least one of the following conditions is satisfied: the problem is determined to be \bfi{unbounded below}, the problem is found to be \bfi{infeasible}, or the convergence ratio $q_s := {(f_s - f_{\mathrm{opt}})}/{(f_0 - f_{\mathrm{opt}})}$ for $s \in \mathcal{S}$ falls below the prescribed tolerance \(\varepsilon = 10^{-4}\) (i.e., \bfi{accuracy reached}). Here, \(\mathcal{S}\) denotes the set of solvers under comparison, \(f_s\) is the best objective value obtained by solver \(s\), \(f_0\) is the objective value at the common starting point, and \(f_{\mathrm{opt}}\) is the best known objective value.
		
		The value \(f_{\mathrm{opt}}\) is typically associated with a global minimizer or, when a global solution is unavailable, with a high-quality local minimizer obtained by applying a collection of local and global derivative-free solvers. We emphasize that \(q_s\) is used solely for benchmarking purposes, since \(f_{\mathrm{opt}}\) is generally unavailable in practical applications.
		
		Otherwise, the problem is classified as \bfi{unsolved} by solver \(s\). This occurs when the solver fails to satisfy any of the above criteria before reaching either the maximum number of function evaluations, denoted by \({\tt nfmax}\), or the maximum permitted computational time, denoted by \({\tt secmax}\). In the present comparison, we set ${\tt nfmax} = 50000$ and ${\tt secmax} = \infty$.
		
		\bfi{Cost measures.}	We use two cost measures: the number of function evaluations,
		\({\tt nf}\), and the elapsed computational time in seconds,
		\({\tt sec}\). Let \(\mathcal{P}\) denote the collection of test problems. For each solver \(s\in\mathcal{S}\) and problem \(p\in\mathcal{P}\), let
		\(c_{s,p}\) denote the computational cost incurred by \(s\) on \(p\).

		\bfi{Solver efficiency and anomaly counters.} For the tabulated comparisons, the efficiency of solver
		\(s\in\mathcal{S}\) on problem \(p\in\mathcal{P}\), with respect to
		\(c_{s,p}\), is defined as
		\[
		e_{s,p}
		:=
		\cases{
			\displaystyle
			100\,
			\frac{\min_{\bar{s}\in\mathcal{S}}c_{\bar{s},p}}
			{c_{s,p}},
			& \mbox{if solver \(s\) solves problem \(p\)},\cr
			0,
			& \mbox{otherwise}.
		}
		\]
		This score measures the performance of \(s\) relative to an ideal
		reference solver that achieves the lowest observed cost on each problem. A value of \(100\%\) indicates the best performance for the selected cost measure; smaller positive values indicate higher relative
		costs, and zero corresponds to an unsuccessful run. The numerical results are summarized in Tables \ref{t.small}--\ref{t.all}. The efficiencies are reported as percentages and are computed separately for
		\(c_{s,p}\in\{{\tt nf},{\tt sec}\}\). In the tables, the corresponding
		percentage efficiency is denoted by \({\tt eff}\%\) for each choice of
		\(c_{s,p}\). Larger values indicate better computational performance. All displayed values are rounded to integers, with nonzero values rounded away from zero when necessary. We also report three anomaly counters,
		$\#n$, $\#t$, and $\#f$.
		These count failures assigned to the three categories
		$n$, $t$, and $f$ by these tables.
		The categories are treated as mutually exclusive.
		Consequently, the number of solved problems plus the
		three anomaly counts equals the number of test problems. Runs terminated because of numerical errors, implementation failures,
		or other algorithm-specific stopping conditions may be reported
		separately.

		\bfi{Solver selection.}
		Table~\ref{a.sol} in Appendix~\ref{sec.inner-solvers} lists a broad
		collection of DFO solvers. Comparing all of them across every test set
		and dimension, however, would make the numerical study unnecessarily
		large and computationally expensive. Therefore, the detailed tables
		focus on a representative subset of solvers selected according to two
		criteria: numerical performance, in terms of robustness and efficiency
		observed in preliminary tests, and methodological diversity. In
		particular, the selected solvers include representatives of model-based,
		line-search, direct-search, and evolution-strategy-type methods.
		
		Nonlinearly constrained solvers are also relevant because linear
		constraints form a special case of nonlinear constraints. We therefore
		include \texttt{COBYLA} \cite{COBYLA} as a well-established benchmark,
		rather than as a method satisfying the feasible-only objective-evaluation
		requirement of \texttt{VLDFO}. Indeed, \texttt{COBYLA} may evaluate
		points that violate the constraints as part of its objective--constraint
		tradeoff. Its inclusion is nevertheless informative because, on the
		considered test problems, it is competitive and in several cases
		outperforms methods designed specifically for linear constraints; see,
		for example, Table~\ref{t.small}. Thus, the comparison with \texttt{COBYLA} provides a performance reference and illustrates the different feasibility requirements of the methods.
		
		Accordingly, the absence of a solver from a particular detailed table
		should not be interpreted as evidence of poor efficiency or robustness;
		rather, it reflects the need to keep the comparisons focused and
		computationally manageable. For the larger-dimensional problems, some
		additional solvers are excluded because their computational cost or
		modelling requirements make them less suitable for the considered
		dimension range.

		\subsection{Infeasibility detection and certificate validation}
		\label{s.infeas-validation}
		
		An important feature of Phase~I of \texttt{VLDFO} is that a problem is
		not declared infeasible merely because a Phase~I subproblem fails to
		produce a feasible point. Instead, \texttt{VLDFO} attempts to obtain and
		validate a dual Farkas-type certificate. The numerical certificate
		validation procedure is described in
		Section~\ref{subsec.vldfo.phase1}. In particular, the direct COI test
		uses the scale-aware conditions \gzit{e.vldfo.coi}, while the
		normalized auxiliary test is based on the quantity $f_I$ defined in
		\gzit{e.getI}, with the scaling quantities and component definitions
		given in \gzit{e.getIscales}. The normalized diagnostic is regarded as successful according to
		\gzit{e.getItest}.
		
		Thus, \texttt{getI} does not replace the direct COI mechanism;
		rather, it provides an additional normalized check of an available
		certificate. A failed \texttt{getI} test is inconclusive and is not
		interpreted as evidence of feasibility. 
		
		The complete experiments produced the Phase~I infeasibility
		classification for the 38 problems listed in
		Table~\ref{t.infeasible-validation}. The same classification is
		reported by all 17 \texttt{VLDFO} variants because Phase~I is common
		to these variants and is performed before the unconstrained
		derivative-free inner solver is invoked. For every listed problem, the stored normalized diagnostic was successful. Thus, these data show the consistency of the common \texttt{VLDFO} feasibility mechanism across the subsequent inner solver configurations; they should not be interpreted as 17 independent infeasibility certificates.
		
		\begin{table}[!ht]
			\caption{Problems classified as infeasible by Phase~I of \texttt{VLDFO}.}
			\label{t.infeasible-validation}
			\centering
			\scriptsize
			\setlength{\tabcolsep}{7pt}
			\begin{tabular}{|lrc|lrc|lrc|}
				\hline
				problem & $n$ & $f_I$ &
				problem & $n$ & $f_I$ &
				problem & $n$ & $f_I$\\
				\hline
				
				\texttt{HS35MOD} & 3 & $0.0000$ &
				\texttt{HS105} & 8 & $6.7477{\times}10^{-19}$ &
				\texttt{EQC} & 9 & $6.1482{\times}10^{-12}$\\
				
				\texttt{QCNEW} & 9 & $6.1482{\times}10^{-12}$ &
				\texttt{CVXQP1\_10} & 10 & $2.0690{\times}10^{-15}$ &
				\texttt{CVXQP3\_10} & 10 & $2.7172{\times}10^{-15}$\\
				
				\texttt{HS112} & 10 & $2.5170{\times}10^{-11}$ &
				\texttt{NCVXQP1\_10} & 10 & $2.0690{\times}10^{-15}$ &
				\texttt{NCVXQP2\_10} & 10 & $2.0690{\times}10^{-15}$\\
				
				\texttt{NCVXQP3\_10} & 10 & $2.0690{\times}10^{-15}$ &
				\texttt{NCVXQP4\_10} & 10 & $2.2281{\times}10^{-15}$ &
				\texttt{NCVXQP6\_10} & 10 & $2.2281{\times}10^{-15}$\\
				
				\texttt{NCVXQP7\_10} & 10 & $2.7172{\times}10^{-15}$ &
				\texttt{NCVXQP8\_10} & 10 & $2.7172{\times}10^{-15}$ &
				\texttt{NCVXQP9\_10} & 10 & $2.7172{\times}10^{-15}$\\
				
				\texttt{GOULDQP1} & 32 & $1.6324{\times}10^{-8}$ &
				\texttt{HIMMELBJ} & 45 & $2.1448{\times}10^{-14}$ &
				\texttt{CVXQP1\_50} & 50 & $8.1614{\times}10^{-16}$\\
				
				\texttt{CVXQP3\_50} & 50 & $8.0944{\times}10^{-16}$ &
				\texttt{NCVXQP1\_50} & 50 & $8.1614{\times}10^{-16}$ &
				\texttt{NCVXQP2\_50} & 50 & $8.1614{\times}10^{-16}$\\
				
				\texttt{NCVXQP3\_50} & 50 & $8.1614{\times}10^{-16}$ &
				\texttt{NCVXQP4\_50} & 50 & $1.6917{\times}10^{-15}$ &
				\texttt{NCVXQP6\_50} & 50 & $1.6917{\times}10^{-15}$\\
				
				\texttt{NCVXQP7\_50} & 50 & $8.0944{\times}10^{-16}$ &
				\texttt{NCVXQP8\_50} & 50 & $8.0944{\times}10^{-16}$ &
				\texttt{NCVXQP9\_50} & 50 & $8.0944{\times}10^{-16}$\\
				
				\texttt{CVXQP1\_100} & 100 & $5.3735{\times}10^{-16}$ &
				\texttt{CVXQP3\_100} & 100 & $1.0607{\times}10^{-15}$ &
				\texttt{NCVXQP1\_100} & 100 & $5.3735{\times}10^{-16}$\\
				
				\texttt{NCVXQP2\_100} & 100 & $5.3735{\times}10^{-16}$ &
				\texttt{NCVXQP3\_100} & 100 & $5.3735{\times}10^{-16}$ &
				\texttt{NCVXQP4\_100} & 100 & $1.4740{\times}10^{-15}$\\
				
				\texttt{NCVXQP6\_100} & 100 & $1.4740{\times}10^{-15}$ &
				\texttt{NCVXQP7\_100} & 100 & $1.0607{\times}10^{-15}$ &
				\texttt{NCVXQP8\_100} & 100 & $1.0607{\times}10^{-15}$\\
				
				\texttt{NCVXQP9\_100} & 100 & $1.0607{\times}10^{-15}$ &
				\texttt{QPCSTAIR} & 467 & $2.8987{\times}10^{-8}$ &
				& &\\
				
				\hline
			\end{tabular}
		\end{table}
		
		Table~\ref{t.infeasible-validation} shows that the infeasibility mechanism
		is exercised over a broad range of dimensions, from $n=3$ to $n=467$. The same Phase~I conclusion is obtained by every
		\texttt{VLDFO} variant. This agreement is expected because the
		infeasibility decision is made before the Phase~II derivative-free solver
		is invoked. Consequently, these 38 cases are counted in the numerical study as
		Phase~I infeasibility terminations of the common \texttt{VLDFO}
		feasibility mechanism rather than as failures of the individual
		derivative-free inner solvers.

		\subsection{Problems with singleton feasible sets}\label{sec:singleton}
		
		For 51 test problems, Phase~I identifies a validated zero-dimensional
		feasible set. Since the feasible region is convex, each such
		zero-dimensional nonempty feasible set consists of a single point.
		Table~\ref{t.singleton} lists these problems together with their
		dimensions.
		
		\begin{center}
			\begin{longtable}{? l@{\hspace{2pt}}r ? l@{\hspace{2pt}}r ? l@{\hspace{2pt}}r ? l@{\hspace{2pt}}r ? l@{\hspace{2pt}}r ?}
				\caption{The 51 test problems for which Phase~I identifies a
					singleton feasible set.}
				\label{t.singleton}\\
				\hline
				problem & $n$ & problem & $n$ & problem & $n$ & problem & $n$ & problem & $n$\\
				\hline
				\endfirsthead
				
				\hline
				\multicolumn{10}{?c?}{Table~\ref{t.singleton} (continued)}\\
				\hline
				problem & $n$ & problem & $n$ & problem & $n$ & problem & $n$ & problem & $n$\\
				\hline
				\endhead
				
				\hline
				\endfoot
				
				\hline
				\endlastfoot
				
				{\tt BOOTH} & 2 &
				{\tt HS24} & 2 &
				{\tt SIPOW1} & 2 &
				{\tt SIPOW1M} & 2 &
				{\tt SIPOW2} & 2\\
				
				{\tt SIPOW2M} & 2 &
				{\tt SUPERSIM} & 2 &
				{\tt TAME} & 2 &
				{\tt ZECEVIC2} & 2 &
				{\tt HS35} & 3\\
				
				{\tt HS35I} & 3 &
				{\tt HS36} & 3 &
				{\tt HS37} & 3 &
				{\tt HS62} & 3 &
				{\tt STANCMIN} & 3\\
				
				{\tt HS44} & 4 &
				{\tt HS44NEW} & 4 &
				{\tt HS76} & 4 &
				{\tt HS76I} & 4 &
				{\tt LIN} & 4\\
				
				{\tt HS86} & 5 &
				{\tt HS55} & 6 &
				{\tt DUALC2} & 7 &
				{\tt AVGASA} & 8 &
				{\tt AVGASB} & 8\\
				
				{\tt DUALC5} & 8 &
				{\tt DUALC8} & 8 &
				{\tt DUALC1} & 9 &
				{\tt HUESTIS\_10} & 10 &
				{\tt LOTSCHD} & 12\\
				
				{\tt PORTFL1} & 12 &
				{\tt PORTFL2} & 12 &
				{\tt PORTFL3} & 12 &
				{\tt PORTFL4} & 12 &
				{\tt PORTFL6} & 12\\
				
				{\tt HS119} & 16 &
				{\tt DEGENLPA} & 20 &
				{\tt DEGENLPB} & 20 &
				{\tt BDRY2\_25} & 25 &
				{\tt LOADBAL} & 31\\
				
				{\tt NASH} & 72 &
				{\tt DUAL4} & 75 &
				{\tt DUAL1} & 85 &
				{\tt DUAL2} & 96 &
				{\tt HIMMELBI} & 100\\
				
				{\tt HUESTIS\_100} & 100 &
				{\tt DUAL3} & 111 &
				{\tt HYDROELS} & 169 &
				{\tt GMNCASE2} & 175 &
				{\tt GMNCASE3} & 175\\
				
				{\tt GMNCASE4} & 175 &
				& &
				& &
				& &
				&\\
				
			\end{longtable}
		\end{center}
		
		\subsection{Results for the dimensions $n\in[1,50]$}\label{s.num1}
		
		The results in Table~\ref{t.small} for dimensions $n\in[1,50]$ show that all 142 test problems are solved by at least one of the tested methods. The highest solved count is obtained by \texttt{VLDFO-CMAES}, which solves all 142 problems and produces no anomalies of any type. It is followed by \texttt{VLDFO-MATRS}, which solves 139 problems; its three remaining cases are due exclusively to reaching the function-evaluation limit. The next highest solved counts are obtained by \texttt{VLDFO-NOMAD-3.9}, with 136 solved problems, and by \texttt{VLDFO-VRBBO} and \texttt{VLDFO-NOMAD-4.6}, both with 134. The two \texttt{NOMAD} variants, however, produce six and eight failure anomalies, respectively, whereas the eight unsolved cases of \texttt{VLDFO-VRBBO} are caused only by reaching the function-evaluation limit.
		
		The main findings from Table~\ref{t.small} can be summarized as follows.
		\begin{itemize}
			\item \texttt{VLDFO-CMAES} solves all problems on this subset. Its coverage exceeds that of the second-best method, \texttt{VLDFO-MATRS}, by about $2.1\%$ of the test set and that of the third-best method, \texttt{VLDFO-NOMAD-3.9}, by about $4.2\%$.
			
			\item Robustness is not confined to a single inner solver. Twelve \texttt{VLDFO} variants exhibit neither failure nor time-limit anomalies; whenever these methods do not solve a problem, the cause is solely exhaustion of the prescribed function-evaluation budget. Thus, on this test set, the framework works with several different inner DFO algorithms.
			
			\item Robustness and evaluation efficiency do not coincide. While \texttt{VLDFO-CMAES} gives the best solved-problem coverage, \texttt{VLDFO-SID-PSM} gives the best function-evaluation efficiency. Its score is two percentage points above the second-best value and three percentage points above the next group of methods, showing that the method solving the largest number of problems is not necessarily the one requiring the fewest evaluations.
			
			\item The fastest methods according to the {\tt sec} efficiency measure are not the most reliable ones. In particular, \texttt{COBYLA} and \texttt{LINCOA} attain high {\tt sec} efficiencies, although their solved-problem coverage is substantially lower and they exhibit considerably more failure anomalies. Their favorable {\tt sec} values may also be influenced by the fact that the underlying solvers are implemented in Fortran and are called from \texttt{MATLAB} through \texttt{MEX} interfaces. Hence, the observed {\tt sec} advantage should not be interpreted solely as an algorithmic advantage, and {\tt sec} efficiency alone gives an incomplete assessment of performance on this test set.
			
			\item Among methods combining high coverage with the absence of failure and time-limit anomalies, \texttt{VLDFO-CMAES}, \texttt{VLDFO-MATRS}, \texttt{VLDFO-VRBBO}, and \texttt{VLDFO-SID-PSM} provide different trade-offs between robustness, solved-problem coverage, and evaluation efficiency.
			
			\item The comparison with direct \texttt{DSPFD} illustrates the effect of the proposed reformulation: \texttt{VLDFO-DSPFD} solves 121 problems compared with only 29 for direct \texttt{DSPFD}. This corresponds to an increase of about $64.8$ percentage points of the complete test set, or more than four times as many solved problems.
		\end{itemize}

		\begin{table}[!ht]
			\caption{Results for dimensions $n\in[1,50]$.}

			\label{t.small}
			
			\begin{center}
				\begin{adjustbox}{center, width=\columnwidth-1pt}
					\footnotesize
					\begin{tabular}{? l? l? l? l? l? l? l?}
						\hline
						\multicolumn{7}{?l?}{stopping test:~~~~~~} \\
						\multicolumn{7}{?c?}{$q_s\le 0.0001,~ {\tt sec}\le \infty,~ {\tt nf}\le 50000$} \\
						\hline
						\multicolumn{7}{?l?}{142 of 142 problems solved} \\
						\cline{1-7}
						\multicolumn{2}{?l?}{dim$\in$[1,50]}&\multicolumn{3}{r?}{\# of anomalies} &\multicolumn{2}{r?}{{\tt eff}\%}\\
						\hline
						\multicolumn{1}{?l?}{solver} & solved&\#n&\#t&\#f&{\tt nf}&{\tt sec}\\
						\hline
						{\tt VLDFO-CMAES}&\cellcolor[HTML]{C0C0C0}142&0&0&0&52&27\\
						{\tt VLDFO-MATRS}&139&3&0&0&55&25\\
						{\tt VLDFO-NOMAD-3.9}&136&0&0&6&56&26\\
						{\tt VLDFO-VRBBO}&134&8&0&0&55&27\\
						{\tt VLDFO-NOMAD-4.6}&134&0&0&8&56&25\\
						{\tt VLDFO-SID-PSM}&133&9&0&0&\cellcolor[HTML]{C0C0C0}59&28\\
						{\tt VLDFO-DFNDFL}&133&9&0&0&54&28\\
						{\tt VLDFO-VRDFON}&132&10&0&0&55&26\\
						{\tt VLDFO-SSDFO}&132&10&0&0&56&28\\
						{\tt VLDFO-SDBOX}&132&10&0&0&56&28\\
						{\tt VLDFO-MADFO}&130&12&0&0&52&27\\
						{\tt VLDFO-ADSMAX}&130&0&0&12&51&27\\
						{\tt VLDFO-MCS}&129&13&0&0&54&25\\
						{\tt VLDFO-GCES}&129&0&0&13&56&27\\
						{\tt COBYLA}&129&3&0&10&41&\cellcolor[HTML]{C0C0C0}64\\
						{\tt VLDFO-NELDER}&128&0&0&14&55&28\\
						{\tt VLDFO-DMADFO}&127&15&0&0&52&27\\
						{\tt VLDFO-FMAES}&125&2&0&15&49&27\\
						{\tt VLDFO-fmincon-sqp}&124&0&0&18&55&34\\
						{\tt VLDFO-UOBYQA}&124&0&0&18&57&27\\
						{\tt VLDFO-fminunc}&122&0&0&20&56&35\\
						{\tt VLDFO-fmincon-active}&122&0&0&20&56&34\\
						{\tt VLDFO-NEWUOA}&122&0&0&20&57&27\\
						{\tt VLDFO-DSPFD}&121&0&0&21&53&27\\
						{\tt VLDFO-NMSMAX}&121&21&0&0&54&27\\
						{\tt LINCOA}&116&1&0&25&41&61\\
						{\tt VLDFO-DFOTR}&115&18&0&9&56&26\\
						{\tt VLDFO-BCDFO}&114&27&0&1&54&26\\
						{\tt DSPFD}&29&3&0&110&15&16\\
						\hline
				\end{tabular}\end{adjustbox}\medskip\end{center}

		\end{table}

		\begin{figure}[!ht]
			\centering
			
			\begin{tikzpicture}
				\begin{axis}[
					ybar,
					bar width=3pt,
					width=\textwidth,
					height=6.5cm,
					title={The five fastest omitted solvers for $n\in[20,30]$},
					ymode=log,
					log origin=infty,
					ymin=0.005,
					ymax=1000,
					ylabel={{\tt sec}},
					ytick={0.01,0.1,1,10,100,1000},
					yticklabels={
						$10^{-2}$,$10^{-1}$,$10^{0}$,
						$10^{1}$,$10^{2}$,$10^{3}$
					},
					xtick={1,...,18},
					xticklabels={
						DTOC3,
						HAGER2,
						DEGENLPA,
						DEGENLPB,
						BDRY225,
						PDE126,
						HAGER4,
						MAKELA4,
						BLOCKQP230,
						BLOCKQP130,
						BLOCKQP530,
						BLOCKQP330,
						BLOWEYA22,
						BLOCKQP430,
						BLOWEYC22,
						BLOWEYB22,
						DEGENQP,
						KSIP
					},
					x tick label style={
						rotate=90,
						anchor=east,
						font=\scriptsize
					},
					tick align=outside,
					grid=major,
					grid style={dashed,gray!25},
					legend style={
						at={(0.5,1.15)},
						anchor=south,
						legend columns=5,
						font=\scriptsize,
						draw=none,
						/tikz/every even column/.append style={column sep=5pt}
					},
					enlarge x limits=0.025
					]
					
					\addplot+[
					ybar,
					bar width=3pt,
					bar shift=-8pt,
					fill=red!75,
					draw=red!90!black,
					forget plot
					] coordinates {
						(1,0.010) (2,0.010) (3,0.020) (4,0.020)
						(5,0.020) (6,0.020) (7,0.020) (8,0.030)
						(9,0.030) (10,0.050) (11,0.110) (12,0.130)
						(13,0.290) (14,0.370) (15,0.600) (16,0.910)
						(17,4.320) (18,9.160)
					};
					
					\addplot+[
					only marks,
					mark=*,
					mark size=2.2pt,
					xshift=-8pt,
					draw=red!90!black,
					fill=white,
					mark options={
						solid,
						draw=red!90!black,
						fill=white,
						line width=0.8pt
					},
					forget plot
					] coordinates {
						(1,0.010) (2,0.010) (3,0.020) (4,0.020)
						(5,0.020) (6,0.020) (7,0.020) (8,0.030)
						(9,0.030) (10,0.050) (11,0.110) (12,0.130)
						(13,0.290) (14,0.370) (15,0.600) (16,0.910)
						(17,4.320) (18,9.160)
					};
					
					\addplot+[
					ybar,
					bar width=3pt,
					bar shift=-4pt,
					fill=blue!70,
					draw=blue!90!black,
					forget plot
					] coordinates {
						(1,0.160) (2,0.170) (3,0.110) (4,0.110)
						(5,0.190) (6,0.110) (7,0.200) (8,0.120)
						(9,2.080) (10,0.320) (11,0.220) (12,0.440)
						(13,1.390) (14,4.480) (15,0.410) (16,0.880)
						(17,64.290) (18,0.620)
					};
					
					\addplot+[
					only marks,
					mark=square*,
					mark size=2.1pt,
					xshift=-4pt,
					mark options={
						solid,
						draw=blue!90!black,
						fill=white,
						line width=0.8pt
					},
					forget plot
					] coordinates {
						(1,0.160) (2,0.170) (3,0.110) (4,0.110)
						(5,0.190) (6,0.110) (7,0.200) (8,0.120)
						(9,2.080) (10,0.320) (11,0.220) (12,0.440)
						(13,1.390) (14,4.480) (15,0.410) (16,0.880)
						(17,64.290) (18,0.620)
					};
					
					\addplot+[
					ybar,
					bar width=3pt,
					bar shift=0pt,
					fill=magenta!70,
					draw=magenta!90!black,
					forget plot
					] coordinates {
						(1,0.160) (2,0.170) (3,0.110) (4,0.110)
						(5,0.190) (6,0.110) (7,0.200) (8,0.100)
						(9,56.800) (10,14.800) (11,10.020) (12,13.410)
						(13,0.400) (14,54.630) (15,1.190) (16,0.390)
						(17,109.060) (18,0.610)
					};
					
					\addplot+[
					only marks,
					mark=triangle*,
					mark size=2.4pt,
					xshift=0pt,
					mark options={
						solid,
						draw=magenta!90!black,
						fill=white,
						line width=0.8pt
					},
					forget plot
					] coordinates {
						(1,0.160) (2,0.170) (3,0.110) (4,0.110)
						(5,0.190) (6,0.110) (7,0.200) (8,0.100)
						(9,56.800) (10,14.800) (11,10.020) (12,13.410)
						(13,0.400) (14,54.630) (15,1.190) (16,0.390)
						(17,109.060) (18,0.610)
					};
					
					\addplot+[
					ybar,
					bar width=3pt,
					bar shift=4pt,
					fill=orange!85,
					draw=orange!65!black,
					forget plot
					] coordinates {
						(1,0.160) (2,0.180) (3,0.110) (4,0.120)
						(5,0.200) (6,0.110) (7,0.210) (8,0.120)
						(9,60.580) (10,0.610) (11,0.420) (12,0.460)
						(13,2.940) (14,105.360) (15,7.720) (16,2.240)
						(17,78.920) (18,2.310)
					};
					
					\addplot+[
					only marks,
					mark=diamond*,
					mark size=2.4pt,
					xshift=4pt,
					mark options={
						solid,
						draw=orange!65!black,
						fill=white,
						line width=0.8pt
					},
					forget plot
					] coordinates {
						(1,0.160) (2,0.180) (3,0.110) (4,0.120)
						(5,0.200) (6,0.110) (7,0.210) (8,0.120)
						(9,60.580) (10,0.610) (11,0.420) (12,0.460)
						(13,2.940) (14,105.360) (15,7.720) (16,2.240)
						(17,78.920) (18,2.310)
					};
					
					\addplot+[
					ybar,
					bar width=3pt,
					bar shift=8pt,
					fill=cyan!70,
					draw=cyan!80!black,
					forget plot
					] coordinates {
						(1,0.180) (2,0.180) (3,0.120) (4,0.110)
						(5,0.200) (6,0.110) (7,0.200) (8,0.120)
						(9,1.780) (10,2.910) (11,1.880) (12,2.380)
						(13,0.230) (14,9.200) (15,0.260) (16,0.200)
						(17,339.840) (18,2.870)
					};
					
					\addplot+[
					only marks,
					mark=star,
					mark size=2.8pt,
					xshift=8pt,
					mark options={
						solid,
						draw=cyan!80!black,
						fill=white,
						line width=0.8pt
					},
					forget plot
					] coordinates {
						(1,0.180) (2,0.180) (3,0.120) (4,0.110)
						(5,0.200) (6,0.110) (7,0.200) (8,0.120)
						(9,1.780) (10,2.910) (11,1.880) (12,2.380)
						(13,0.230) (14,9.200) (15,0.260) (16,0.200)
						(17,339.840) (18,2.870)
					};
					
					\addlegendimage{
						red!90!black,
						line width=1.5pt,
						mark=*,
						mark size=2.2pt,
						mark options={draw=red!90!black,fill=white}
					}
					\addlegendentry{\texttt{COBYLA}}
					
					\addlegendimage{
						blue!90!black,
						line width=1.5pt,
						mark=square*,
						mark size=2.1pt,
						mark options={draw=blue!90!black,fill=white}
					}
					\addlegendentry{\texttt{VLDFO-DFOTR}}
					
					\addlegendimage{
						magenta!90!black,
						line width=1.5pt,
						mark=triangle*,
						mark size=2.4pt,
						mark options={draw=magenta!90!black,fill=white}
					}
					\addlegendentry{\texttt{VLDFO-GCES}}
					
					\addlegendimage{
						orange!65!black,
						line width=1.5pt,
						mark=diamond*,
						mark size=2.4pt,
						mark options={draw=orange!65!black,fill=white}
					}
					\addlegendentry{\texttt{VLDFO-BCDFO}}
					
					\addlegendimage{
						cyan!80!black,
						line width=1.5pt,
						mark=star,
						mark size=2.8pt,
						mark options={draw=cyan!80!black,fill=white}
					}
					\addlegendentry{\texttt{VLDFO-UOBYQA}}
					
				\end{axis}
			\end{tikzpicture}
			
			\caption{Per-problem {\tt sec} values for the five fastest methods among the
				12 solvers omitted from the larger-dimensional comparison, restricted
				to the 18 feasible test problems with dimensions $n\in[20,30]$, ordered from left to right by increasing
				{\tt sec} of \texttt{COBYLA}. The ranking is based on mean {\tt sec} over these problems and is
				restricted to methods having positive recorded timing data for all
				18 problems. The five fastest methods are
				\texttt{COBYLA}, \texttt{VLDFO-DFOTR},
				\texttt{VLDFO-GCES}, \texttt{VLDFO-BCDFO}, and
				\texttt{VLDFO-UOBYQA}, with mean {\tt sec} values of approximately
				$0.90$, $4.24$, $14.58$, $14.60$, and $20.15$,
				respectively. Each method is identified by
				both a distinct color and a distinct marker located at the upper
				end of each bar.}
			\label{fig:five-fastest-20-30}
		\end{figure}
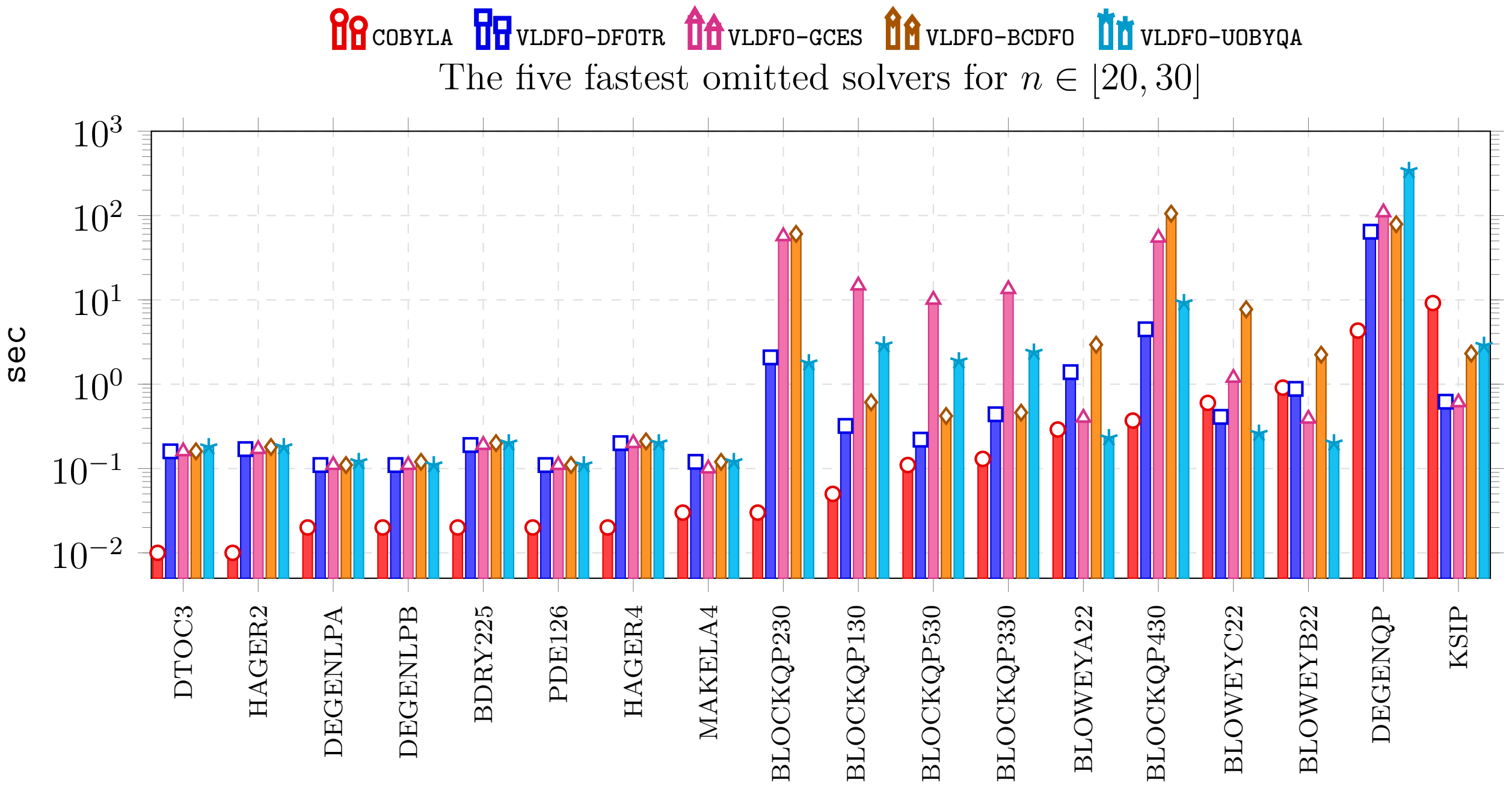

		Figure~\ref{fig:five-fastest-20-30} shows the timing data used in the solver selection for the larger-dimensional experiments. Even among the five fastest methods from the omitted set having positive recorded timing data for all 18 problems, the {\tt sec} values exhibit considerable variability on problems with $n\in[20,30]$, and several individual runs require tens or hundreds of seconds. Together with the preliminary results concerning robustness and scalability, these results motivated the decision to omit these 12 methods from the subsequent experiments with $n\in[51,500]$, while retaining a representative set of methods for the larger-dimensional comparison.

		\subsection{Results for the dimensions $n\in[1,500]$}\label{s.all}
		
		Table~\ref{t.all} reports the results for the complete set of problems with dimensions $n\in[1,500]$ together with, in the right panel, the corresponding results for the larger-dimensional subset $n\in[51,500]$. Over the complete test set, 220 of the 246 problems are solved by at least one method, while the remaining 26 problems are not solved by any method within the prescribed function-evaluation budget. The largest individual solved count is achieved by \texttt{VLDFO-CMAES}, which solves 209 problems without producing any failure or time-limit anomalies. It is followed by \texttt{VLDFO-VRBBO}, with 201 solved problems, and \texttt{VLDFO-VRDFON}, with 200. The next highest solved counts among the variants without failure or time-limit anomalies are obtained by \texttt{VLDFO-DFNDFL}, with 197 solved problems, and \texttt{VLDFO-SDBOX}, with 196.
		
		The main findings from Table~\ref{t.all} can be summarized as follows.
		\begin{itemize}
			\item Overall, about $89.4\%$ of the complete test set is solved by at least one
			method. \texttt{VLDFO-CMAES} alone solves about $85.0\%$ of all
			problems and $95.0\%$ of the 220 problems solved by at least one
			method. It solves eight more problems than the second-best method,
			corresponding to about $3.3$ percentage points of the complete test set.
			
			\item The methods with the highest solved counts also show good numerical reliability. In particular, \texttt{VLDFO-CMAES}, \texttt{VLDFO-VRBBO},
			\texttt{VLDFO-VRDFON}, \texttt{VLDFO-DFNDFL}, and
			\texttt{VLDFO-SDBOX} combine high solved counts with no failure
			or time-limit anomalies.
			
			\item The differences in function-evaluation efficiency among the
			leading methods are small. \texttt{VLDFO-SDBOX} and
			\texttt{VLDFO-SSDFO} attain the highest {\tt nf} efficiency,
			only one percentage point above \texttt{VLDFO-CMAES},
			\texttt{VLDFO-VRBBO}, and \texttt{VLDFO-VRDFON}.
			Thus, several of the most reliable methods are also competitive
			in evaluation cost.
			
			\item The highest {\tt sec} efficiencies are obtained by variants based
			on \texttt{fmincon} and \texttt{fminunc}, but these methods also
			produce many failure anomalies. Among the variants with no failure
			or time-limit anomalies, \texttt{VLDFO-CMAES},
			\texttt{VLDFO-DFNDFL}, and \texttt{VLDFO-SDBOX} provide the
			best timing efficiency. Hence, the results favor methods that
			balance speed and reliability rather than timing alone.
			
			\item The larger-dimensional results in the right panel show a broadly
			similar behavior. Of the 104 problems with $n\in[51,500]$, 78 are
			solved by at least one method. The main difference is that
			\texttt{VLDFO-VRDFON} becomes the leading method, solving 68 problems,
			compared with 67 for both \texttt{VLDFO-VRBBO} and
			\texttt{VLDFO-CMAES}. It also attains the highest {\tt nf} efficiency,
			equal to $49\%$, with no failure or time-limit anomalies. Thus, together
			with Table~\ref{t.small}, the results indicate that
			\texttt{VLDFO-CMAES} performs best overall and on smaller problems,
			whereas \texttt{VLDFO-VRDFON} has a slight advantage on the
			larger-dimensional subset.
		\end{itemize}

		\begin{table}[!ht]
			\caption{Results for dimensions $n\in[1,500]$ (left) and for the
				larger-dimensional subset $n\in[51,500]$ (right).}
			\label{t.all}
			
			\begin{center}
				\begin{minipage}[t]{0.48\textwidth}
					\centering
					\begin{adjustbox}{center, width=\linewidth}
						\scriptsize
						\scalebox{0.92}{\begin{tabular}{? l? l? l? l? l? l? l?}
								\hline
								\multicolumn{7}{?l?}{stopping test:~~~~~~} \\
								\multicolumn{7}{?c?}{$q_s\le 0.0001,~
									{\tt sec}\le \infty,~{\tt nf}\le 50000$} \\
								\hline
								\multicolumn{7}{?l?}{220 of 246 problems solved} \\
								\cline{1-7}
								\multicolumn{2}{?l?}{dim$\in[1,500]$}&
								\multicolumn{3}{r?}{\# of anomalies}&
								\multicolumn{2}{r?}{{\tt eff}\%}\\
								\hline
								solver&solved&\#n&\#t&\#f&{\tt nf}&{\tt sec}\\
								\hline
								{\tt VLDFO-CMAES}&\cellcolor[HTML]{C0C0C0}209&37&0&0&54&63\\
								{\tt VLDFO-VRBBO}&201&45&0&0&54&58\\
								{\tt VLDFO-VRDFON}&200&46&0&0&54&58\\
								{\tt VLDFO-DFNDFL}&197&49&0&0&53&63\\
								{\tt VLDFO-SDBOX}&196&50&0&0&
								\cellcolor[HTML]{C0C0C0}55&63\\
								{\tt VLDFO-ADSMAX}&193&28&0&25&50&59\\
								{\tt VLDFO-SSDFO}&192&54&0&0&
								\cellcolor[HTML]{C0C0C0}55&60\\
								{\tt VLDFO-MADFO}&187&59&0&0&49&57\\
								{\tt VLDFO-MCS}&187&59&0&0&50&52\\
								{\tt VLDFO-FMAES}&186&35&0&25&47&58\\
								{\tt VLDFO-NELDER}&186&3&0&57&51&58\\
								{\tt VLDFO-DMADFO}&184&62&0&0&48&56\\
								{\tt VLDFO-fmincon-sqp}&183&1&0&62&54&
								\cellcolor[HTML]{C0C0C0}70\\
								{\tt VLDFO-fminunc}&179&1&0&66&54&
								\cellcolor[HTML]{C0C0C0}70\\
								{\tt VLDFO-fmincon-active}&178&7&0&61&53&69\\
								{\tt VLDFO-DSPFD}&178&3&0&65&50&57\\
								{\tt VLDFO-NMSMAX}&178&68&0&0&51&58\\
								\hline
						\end{tabular}}
					\end{adjustbox}
				\end{minipage}
				\hfill
				\begin{minipage}[t]{0.48\textwidth}
					\centering
					\begin{adjustbox}{center, width=\linewidth}
						\scriptsize
						\scalebox{0.92}{\begin{tabular}{? l? l? l? l? l? l? l?}
								\hline
								\multicolumn{7}{?l?}{stopping test:~~~~~~} \\
								\multicolumn{7}{?c?}{$q_s\le 0.0001,~
									{\tt sec}\le \infty,~{\tt nf}\le 50000$} \\
								\hline
								\multicolumn{7}{?l?}{78 of 104 problems solved} \\
								\cline{1-7}
								\multicolumn{2}{?l?}{dim$\in[51,500]$}&
								\multicolumn{3}{r?}{\# of anomalies}&
								\multicolumn{2}{r?}{{\tt eff}\%}\\
								\hline
								solver&solved&\#n&\#t&\#f&{\tt nf}&{\tt sec}\\
								\hline
								{\tt VLDFO-VRDFON}&\cellcolor[HTML]{C0C0C0}68&36&0&0&
								\cellcolor[HTML]{C0C0C0}49&52\\
								{\tt VLDFO-VRBBO}&67&37&0&0&45&49\\
								{\tt VLDFO-CMAES}&67&37&0&0&47&
								\cellcolor[HTML]{C0C0C0}53\\
								{\tt VLDFO-SDBOX}&64&40&0&0&47&
								\cellcolor[HTML]{C0C0C0}53\\
								{\tt VLDFO-DFNDFL}&64&40&0&0&45&52\\
								{\tt VLDFO-ADSMAX}&63&28&0&13&42&49\\
								{\tt VLDFO-FMAES}&61&33&0&10&40&47\\
								{\tt VLDFO-SSDFO}&60&44&0&0&47&51\\
								{\tt VLDFO-fmincon-sqp}&59&1&0&44&44&
								\cellcolor[HTML]{C0C0C0}54\\
								{\tt VLDFO-MCS}&58&46&0&0&40&42\\
								{\tt VLDFO-NELDER}&58&3&0&43&38&46\\
								{\tt VLDFO-MADFO}&57&47&0&0&40&46\\
								{\tt VLDFO-fminunc}&57&1&0&46&43&53\\
								{\tt VLDFO-DSPFD}&57&3&0&44&41&47\\
								{\tt VLDFO-NMSMAX}&57&47&0&0&41&47\\
								{\tt VLDFO-DMADFO}&57&47&0&0&40&46\\
								{\tt VLDFO-fmincon-active}&56&7&0&41&43&53\\
								\hline
						\end{tabular}}
					\end{adjustbox}
				\end{minipage}
				\medskip
			\end{center}
		\end{table}
		
		\bfi{Solved-case diagnostics.}
		Table~\ref{t.infeas} provides a more detailed breakdown of the solved cases for dimensions $n\in[1,500]$, distinguishing solutions obtained by satisfying the target-accuracy criterion, cases classified as unbounded below, and problems identified as infeasible by the common \texttt{VLDFO} feasibility mechanism.

		\begin{table}[!ht]
			\caption{Solved-case diagnostics for dimensions $n\in[1,500]$,
				showing the contributions of target-accuracy attainment,
				the unbounded-below threshold, and infeasibility detection
				to the reported solved counts.}
			\label{t.infeas}
			
			\begin{center}
				\begin{adjustbox}{center, width=\columnwidth-1pt}
					\footnotesize
					\begin{tabular}{? l? r? r? r? r?}
						\hline
						\multicolumn{5}{?c?}{Solved-case diagnostics: 220 of 246 problems solved}\\
						\hline
						\multicolumn{5}{?c?}{dim$\in[1,500]$}\\
						\hline
						solver& solved& $q_s \leq$ tol& $f_{\rm best}\leq -10^{12}$& err = infeasible\\
						\hline
						{\tt VLDFO-CMAES}&209&170&1&38\\
						{\tt VLDFO-VRBBO}&201&162&1&38\\
						{\tt VLDFO-VRDFON}&200&161&1&38\\
						{\tt VLDFO-DFNDFL}&197&158&1&38\\
						{\tt VLDFO-SDBOX}&196&157&1&38\\
						{\tt VLDFO-ADSMAX}&193&154&1&38\\
						{\tt VLDFO-SSDFO}&192&154&0&38\\
						{\tt VLDFO-MADFO}&187&149&0&38\\
						{\tt VLDFO-MCS}&187&148&1&38\\
						{\tt VLDFO-FMAES}&186&148&0&38\\
						{\tt VLDFO-NELDER}&186&148&0&38\\
						{\tt VLDFO-DMADFO}&184&145&1&38\\
						{\tt VLDFO-fmincon-sqp}&183&144&1&38\\
						{\tt VLDFO-fminunc}&179&141&0&38\\
						{\tt VLDFO-fmincon-active}&178&140&0&38\\
						{\tt VLDFO-DSPFD}&178&140&0&38\\
						{\tt VLDFO-NMSMAX}&178&140&0&38\\
						\hline
					\end{tabular}
				\end{adjustbox}
				\medskip
			\end{center}
		\end{table}
		
		\bfi{List of unsolved problems.}
		Table~\ref{t.unsolved} lists all 26 test problems
		that remain unsolved by every tested solver for
		dimensions $n\in[1,500]$. Of the 246 test problems,
		220 are solved by at least one method, while the
		remaining 26 do not satisfy the prescribed
		solved-case criteria within the computational limits.

		\clearpage
		
		\begin{center}
			\begin{longtable}{?l l l l?}
				\caption{List of the 26 test problems unsolved by all
					tested solvers for dimensions $n\in[1,500]$.}
				\label{t.unsolved}\\
				\hline
				\multicolumn{4}{?c?}{26 test problems unsolved by all solvers, $n\in[1,500]$}\\
				\hline
				\endfirsthead
				
				\hline
				\multicolumn{4}{?c?}{Table~\ref{t.unsolved} (continued)}\\
				\hline
				\endhead
				
				\hline
				\multicolumn{4}{??c??}{Continued on next page}\\
				\hline
				\endfoot
				
				\hline
				\endlastfoot
				
				{\tt STNQP165}&{\tt GOULDQP299}&{\tt LISWET9}&{\tt LISWET11102}\\
				{\tt LISWET1102}&{\tt LISWET2102}&{\tt LISWET4102}&{\tt LISWET5102}\\
				{\tt LISWET7102}&{\tt READING2153}&{\tt GOULDQP2199}&{\tt GOULDQP3199}\\
				{\tt PRIMALC5}&{\tt GOULDQP2299}&{\tt PRIMAL1}&{\tt GOULDQP2399}\\
				{\tt LISWET10402}&{\tt LISWET11402}&{\tt LISWET12402}&{\tt LISWET1402}\\
				{\tt LISWET2402}&{\tt LISWET4402}&{\tt LISWET6402}&{\tt LISWET8402}\\
				{\tt LISWET9402}&{\tt GOULDQP2499}&&\\
				
			\end{longtable}
		\end{center}

		\subsection{Comparison of the \texttt{NOMAD} variants for the dimensions $n\in[1,10]$}
		
		To examine more closely the role of \texttt{NOMAD} inside the proposed
		framework, Table~\ref{t.nomad} compares three configurations on
		the subset of 85 problems with dimensions $n\in[1,10]$.
		The first two methods, \texttt{VLDFO-NOMAD-4.6} and
		\texttt{VLDFO-NOMAD-3.9}, use \texttt{NOMAD} only as the unconstrained
		inner solver of \texttt{VLDFO}. In contrast, \texttt{NOMAD-4.6}
		is applied directly to the original linearly constrained
		optimization problem.
		
		The difference is substantial. Both \texttt{VLDFO-NOMAD-4.6}
		and \texttt{VLDFO-NOMAD-3.9} solve all 85 problems without any
		anomaly, whereas the direct linearly constrained
		\texttt{NOMAD-4.6} solves only 68 problems and produces
		17 failure anomalies. The two \texttt{VLDFO} variants are also
		considerably more efficient in terms of {\tt nf}: their corresponding
		{\tt nf} efficiency scores are 77\% and 76\%, compared with only
		28\% for the direct linearly constrained version.
		For {\tt sec}, \texttt{VLDFO-NOMAD-3.9} gives the best efficiency
		score, equal to 71\%, followed by \texttt{VLDFO-NOMAD-4.6}
		with 57\%, whereas the direct \texttt{NOMAD-4.6} obtains 44\%.
		
		\begin{table}[!ht]
			\caption{Comparison of two \texttt{NOMAD} versions used as unconstrained
				inner solvers of \texttt{VLDFO} with the direct linearly
				constrained \texttt{NOMAD-4.6} solver for dimensions
				$n\in[1,10]$.}
			\label{t.nomad}
			
			\begin{center}
				\begin{adjustbox}{center, width=\columnwidth-1pt}
					\footnotesize
					\begin{tabular}{? l? l? l? l? l? l? l? l? l?}
						\hline
						\multicolumn{7}{?l?}{stopping test:~~~~~~} \\
						\multicolumn{7}{?c?}{$q_s\le 0.0001,~ {\tt sec}\le \infty,~
							{\tt nf}\le 50000$} \\
						\hline
						\multicolumn{7}{?l?}{85 test problems considered} \\
						\cline{1-7}
						\multicolumn{2}{?l?}{dim$\in[1,10]$}&
						\multicolumn{3}{r?}{\# of anomalies} &
						\multicolumn{2}{r?}{{\tt eff}\%}\\
						\hline
						\multicolumn{1}{?l?}{solver} & solved&\#n&\#t&\#f&
						{\tt nf}&{\tt sec}\\
						\hline
						\texttt{VLDFO-NOMAD-4.6}&
						\cellcolor[HTML]{C0C0C0}85&0&0&0&
						\cellcolor[HTML]{C0C0C0}77&57\\
						\texttt{VLDFO-NOMAD-3.9}&
						\cellcolor[HTML]{C0C0C0}85&0&0&0&76&
						\cellcolor[HTML]{C0C0C0}71\\
						\texttt{NOMAD-4.6}&68&0&0&17&28&44\\
						\hline
					\end{tabular}
				\end{adjustbox}
				\medskip
			\end{center}
		\end{table}
		
\begin{figure}[!ht]
			\centering
			
			\begin{tikzpicture}
				\begin{axis}[
					ybar,
					bar width=2pt,
					width=\textwidth,
					height=4cm,
					extra description/.code={
						\node[font=\small,anchor=south]
						at (axis description cs:0.5,1.08)
						{Problems with $n\in[2,4]$};
					},
					ymode=log,
					log origin=infty,
					ymin=0.001,
					ymax=10000,
					ylabel={{\tt sec}},
					ytick={0.001,0.01,0.1,1,10,100,1000,10000},
					yticklabels={$10^{-3}$,$10^{-2}$,$10^{-1}$,$10^{0}$,$10^{1}$,$10^{2}$,$10^{3}$,$10^{4}$},
					xtick={1,...,37},
					xticklabels={
						EXTRASIM,OET1,OET3,HUBFIT,LSQFIT,PT,SUPERSIM,TAME,
						ZECEVIC2,TFI2,HS44,SIMPLLPA,SIMPLLPB,HS35I,STANCMIN,
						HIMMELBA,HS21,HS24,HS36,HS62,HS37,HS35,HS41,HS76,
						HS76I,BOOTH,HS44NEW,LIN,ZANGWIL3,HATFLDH,BIGGSC4,
						SIPOW1M,SIPOW2,SIPOW2M,SIPOW1,TFI3,SIPOW4
					},
					x tick label style={rotate=90,anchor=east,font=\tiny},
					tick align=outside,
					grid=major,
					grid style={dashed,gray!25},
					legend style={
						at={(0.5,1.30)},
						anchor=south,
						legend columns=3,
						font=\scriptsize,
						draw=none
					},
					enlarge x limits=0.01
					]
					
					\addplot+[
					ybar,
					bar width=2pt,
					bar shift=-2.6pt
					] coordinates {
						(1,0.340) (2,1.570) (3,1.450) (4,0.690) (5,0.660)
						(6,0.590) (7,0.020) (8,0.100) (9,0.100) (10,0.800)
						(11,0.070) (12,0.560) (13,0.640) (14,0.090) (15,0.130)
						(16,0.040) (17,0.650) (18,0.130) (19,0.070) (20,0.080)
						(21,0.070) (22,0.080) (23,0.780) (24,0.070) (25,0.200)
						(26,0.030) (27,0.090) (28,0.190) (29,0.050) (30,2.900)
						(31,0.820) (32,3.670) (33,2.010) (34,2.640) (35,3.630)
						(36,0.690) (37,8.220)
					};
					
					\addplot+[
					ybar,
					bar width=2pt,
					bar shift=0pt
					] coordinates {
						(1,0.330) (2,1.170) (3,1.220) (4,0.450) (5,0.260)
						(6,0.360) (7,0.030) (8,0.100) (9,0.110) (10,0.350)
						(11,0.070) (12,0.220) (13,0.220) (14,0.100) (15,0.130)
						(16,0.040) (17,0.350) (18,0.130) (19,0.070) (20,0.070)
						(21,0.080) (22,0.090) (23,0.490) (24,0.060) (25,0.210)
						(26,0.020) (27,0.090) (28,0.180) (29,0.040) (30,1.220)
						(31,0.560) (32,3.660) (33,2.020) (34,2.640) (35,3.620)
						(36,0.320) (37,2.700)
					};
					
					\addplot+[
					ybar,
					bar width=2pt,
					bar shift=2.6pt
					] coordinates {
						(1,0.005) (2,0.005) (3,0.005) (4,0.010) (5,0.010)
						(6,0.010) (7,0.010) (8,0.010) (9,0.010) (10,0.010)
						(11,0.010) (12,0.020) (13,0.020) (14,0.020) (15,0.020)
						(16,0.030) (17,0.030) (18,0.030) (19,0.050) (20,0.080)
						(21,0.090) (22,0.110) (23,0.130) (24,0.140) (25,0.150)
						(26,0.160) (27,0.210) (28,0.490) (29,0.560) (30,2.180)
						(31,4.250) (32,4.430) (33,4.430) (34,4.870) (35,5.000)
						(36,9.230) (37,336.220)
					};
					
					\legend{
						\texttt{VLDFO-NOMAD-4.6},
						\texttt{VLDFO-NOMAD-3.9},
						\texttt{NOMAD-4.6}
					}
				\end{axis}
			\end{tikzpicture}
			
			\vspace{0.4cm}
			
			\begin{minipage}[t]{0.49\textwidth}
				\centering
				
				\begin{tikzpicture}
					\begin{axis}[
						ybar,
						bar width=2pt,
						width=\linewidth,
						height=4cm,
						title={Problems with $n\in[5,7]$},
						ymode=log,
						log origin=infty,
						ymin=0.001,
						ymax=10000,
						ylabel={{\tt sec}},
						ytick={0.001,0.01,0.1,1,10,100,1000,10000},
						yticklabels={$10^{-3}$,$10^{-2}$,$10^{-1}$,$10^{0}$,$10^{1}$,$10^{2}$,$10^{3}$,$10^{4}$},
						xtick={1,...,19},
						xticklabels={
							DUALC2,HS86,HS54,HS55,HS21MOD,STCQP15,HS48,
							STNQP15,STNQP25,HS52,HS50,HS49,HS53,BT3,
							EXPFITA,HS51,EXPFITB,EXPFITC,HS268
						},
						x tick label style={rotate=90,anchor=east,font=\tiny},
						tick align=outside,
						grid=major,
						grid style={dashed,gray!25},
						enlarge x limits=0.02
						]
						
						\addplot+[
						ybar,
						bar width=2pt,
						bar shift=-2.6pt
						] coordinates {
							(1,0.290) (2,0.110) (3,0.410) (4,0.160)
							(5,0.390) (6,0.410) (7,0.640) (8,0.500)
							(9,0.530) (10,0.720) (11,0.470) (12,0.570)
							(13,0.570) (14,0.830) (15,3.240) (16,0.550)
							(17,1.860) (18,4.420) (19,2.090)
						};
						
						\addplot+[
						ybar,
						bar width=2pt,
						bar shift=0pt
						] coordinates {
							(1,0.290) (2,0.100) (3,0.250) (4,0.170)
							(5,0.250) (6,0.220) (7,0.260) (8,0.300)
							(9,0.300) (10,0.260) (11,0.230) (12,0.250)
							(13,0.230) (14,0.250) (15,1.580) (16,0.280)
							(17,1.420) (18,1.710) (19,1.050)
						};
						
						\addplot+[
						ybar,
						bar width=2pt,
						bar shift=2.6pt
						] coordinates {
							(1,0.010) (2,0.020) (3,0.050) (4,0.220)
							(5,0.510) (6,0.710) (7,1.430) (8,1.540)
							(9,1.610) (10,4.410) (11,4.440) (12,6.270)
							(13,9.540) (14,9.950) (15,11.850) (16,12.130)
							(17,42.310) (18,80.320) (19,4781.350)
						};
						
					\end{axis}
				\end{tikzpicture}
				
			\end{minipage}
			\hfill
			\begin{minipage}[t]{0.49\textwidth}
				\centering
				
				\begin{tikzpicture}
					\begin{axis}[
						ybar,
						bar width=2pt,
						width=\linewidth,
						height=4cm,
						title={Problems with $n\in[8,10]$},
						ymode=log,
						log origin=infty,
						ymin=0.001,
						ymax=10000,
						ylabel={{\tt sec}},
						ytick={0.001,0.01,0.1,1,10,100,1000,10000},
						yticklabels={$10^{-3}$,$10^{-2}$,$10^{-1}$,$10^{0}$,$10^{1}$,$10^{2}$,$10^{3}$,$10^{4}$},
						xtick={1,...,15},
						xticklabels={
							ARGLBLE10,DUALC5,DUALC8,DUALC1,POWELL20,
							HUESTIS10,AVGASB,AVGASA,HS112,READING29,
							ODFITS,GENHS28,ARGLALE10,ARGLCLE10,DEGENQP10
						},
						x tick label style={rotate=90,anchor=east,font=\tiny},
						tick align=outside,
						grid=major,
						grid style={dashed,gray!25},
						enlarge x limits=0.02
						]
						
						\addplot+[
						ybar,
						bar width=2pt,
						bar shift=-2.6pt
						] coordinates {
							(1,1.030) (2,0.220) (3,0.400) (4,0.240)
							(5,0.430) (6,0.090) (7,0.100) (8,0.090)
							(9,0.005) (10,0.530) (11,0.450) (12,0.410)
							(13,0.860) (14,0.920) (15,4.850)
						};
						
						\addplot+[
						ybar,
						bar width=2pt,
						bar shift=0pt
						] coordinates {
							(1,0.850) (2,0.210) (3,0.400) (4,0.250)
							(5,0.270) (6,0.090) (7,0.100) (8,0.100)
							(9,0.005) (10,0.300) (11,0.220) (12,0.190)
							(13,0.880) (14,0.750) (15,5.720)
						};
						
						\addplot+[
						ybar,
						bar width=2pt,
						bar shift=2.6pt
						] coordinates {
							(1,0.005) (2,0.010) (3,0.010) (4,0.010)
							(5,0.010) (6,0.650) (7,1.070) (8,1.190)
							(9,2.230) (10,2.660) (11,7.970) (12,16.060)
							(13,37.100) (14,45.060) (15,128.210)
						};
						
					\end{axis}
				\end{tikzpicture}
				
			\end{minipage}
			
			\caption{Comparison of {\tt sec} for the feasible problems with dimensions $n\in[1,10]$
				for which recorded {\tt sec} data are available for at least one of the three
				solvers. The problems are grouped by dimension ranges: $n\in[2,4]$,
				$n\in[5,7]$, and $n\in[8,10]$. Within each dimension range, the problems
				are ordered from left to right by increasing {\tt sec} of the direct linearly
				constrained \texttt{NOMAD-4.6}, which is used as the reference solver for
				the ordering. For each problem, the three adjacent bars correspond to
				\texttt{VLDFO-NOMAD-4.6}, \texttt{VLDFO-NOMAD-3.9}, and the direct
				linearly constrained \texttt{NOMAD-4.6}, respectively. A logarithmic
				scale is used because the recorded {\tt sec} values range from as low as
				$5\times10^{-3}$ seconds to $4781.35$ seconds.}
			\label{fig:nomad-time}
		\end{figure}
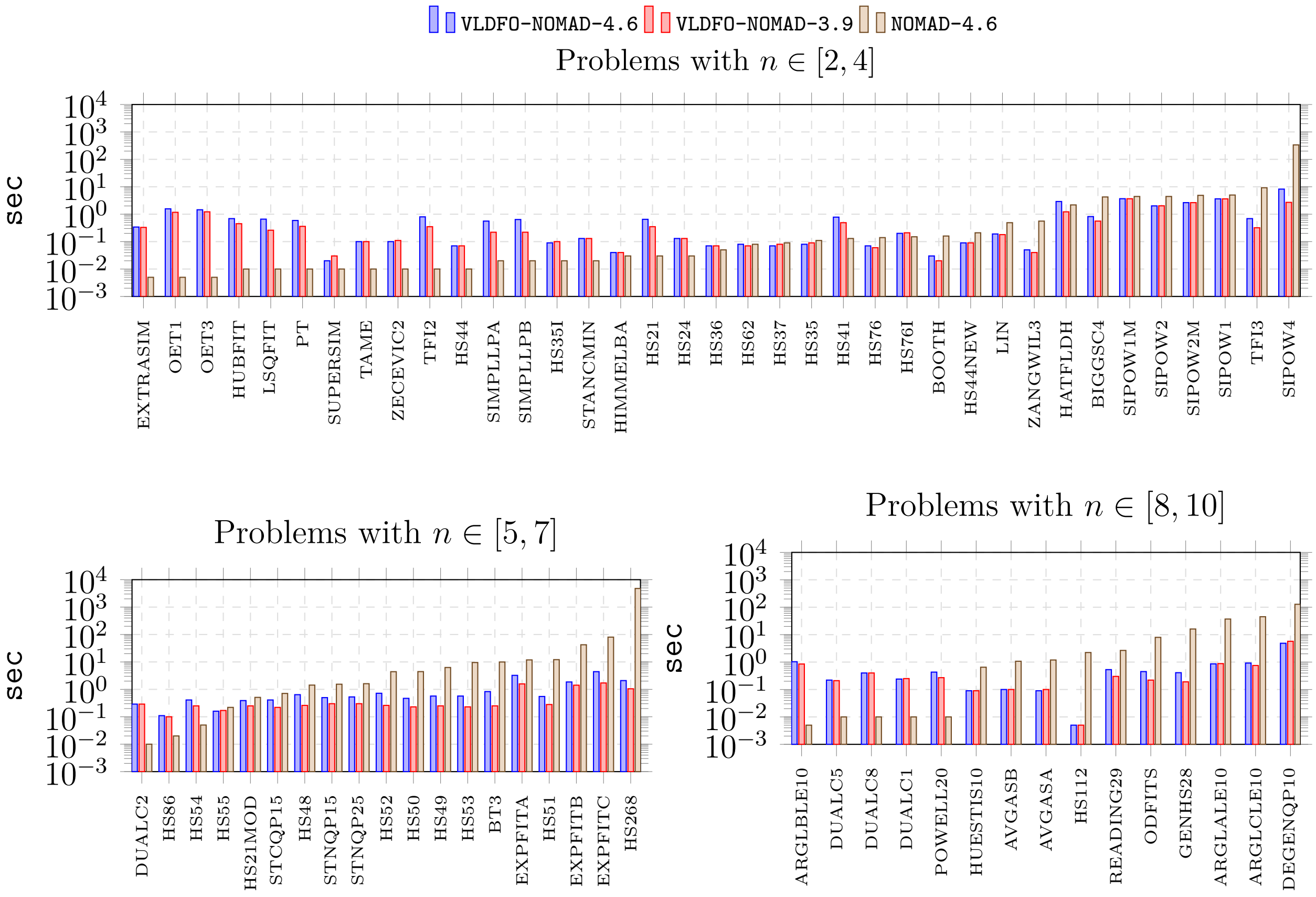

		The {\tt sec} data give an even clearer picture. For the runs with
		finite recorded times, the mean {\tt sec} values are approximately
		$0.757$ for \texttt{VLDFO-NOMAD-4.6},
		$0.512$ for \texttt{VLDFO-NOMAD-3.9}, and
		$11.764$ for the direct linearly constrained
		\texttt{NOMAD-4.6}; see Figure~\ref{fig:nomad-time}.
		Thus, on this comparison, the mean {\tt sec} value of the direct
		linearly constrained \texttt{NOMAD-4.6} is about $15.5$ times that
		of \texttt{VLDFO-NOMAD-4.6} and about $23.0$ times that of
		\texttt{VLDFO-NOMAD-3.9}.
		
		\bfi{Infeasible problems for $n\in[1,10]$.}
		The following 15 problems are identified as infeasible by the common
		\texttt{VLDFO} feasibility mechanism:
		\texttt{HS35MOD}, \texttt{HS105}, \texttt{EQC}, \texttt{QCNEW},
		\texttt{CVXQP110}, \texttt{CVXQP310}, \texttt{HS112},
		\texttt{NCVXQP110}, \texttt{NCVXQP210}, \texttt{NCVXQP310},
		\texttt{NCVXQP410}, \texttt{NCVXQP610}, \texttt{NCVXQP710},
		\texttt{NCVXQP810}, and \texttt{NCVXQP910}. For the direct linearly constrained \texttt{NOMAD-4.6} runs, two
		termination outcomes may occur when only infeasible points are obtained.
		\texttt{RunFlag=-1} indicates that the mesh stopping criterion has been
		reached and that a local minimum of the infeasibility measure has been
		obtained; for smooth constraint functions, this provides only a local
		infeasibility conclusion. In contrast, \texttt{RunFlag=-2} indicates
		that only infeasible points were obtained without reaching the mesh
		criterion and provides no corresponding local-convergence guarantee.
		Thus, neither outcome constitutes a global certificate that the
		original optimization problem is infeasible. Consequently, these
		\texttt{NOMAD-4.6} outcomes are not counted as infeasibility detections
		in the present benchmark.

\section{Conclusion}

We proposed \texttt{VLDFO}, a feasible-point framework for
DFO with linear constraints. Phase~I uses
only the constraint data to detect infeasibility or a singleton feasible set
and, when the feasible set has positive dimension, to identify implied
activities, determine the affine subspace spanned by the feasible set, and
construct a validated reduced representation.
Phase~II then transforms the reduced
linearly constrained problem into an unconstrained problem through a
pointed feasible retraction. Consequently, the original objective
oracle is evaluated only at points that pass the feasibility checks,
while a broad range of unconstrained DFO methods can be
used as inner solvers.

The numerical experiments show that performance depends on the
choice of inner solver. \texttt{VLDFO-CMAES} solves the largest number
of problems over the complete test set, whereas
\texttt{VLDFO-VRDFON} solves the largest number on the
larger-dimensional subset.
\texttt{VLDFO-VRBBO}, \texttt{VLDFO-SDBOX}, and
\texttt{VLDFO-DFNDFL} also achieve competitive solved counts and
efficiency scores. On the present test set, the two
\texttt{VLDFO-NOMAD} variants solve more of the selected problems
than direct linearly constrained \texttt{NOMAD-4.6} and have lower
recorded mean CPU times, as measured by {\tt sec}.

The exact-penalty results concern the exact transformed objective.
Thus, the transformation provides a common feasible-evaluation
framework, while convergence and practical performance depend on
the assumptions and mechanisms of the selected unconstrained DFO
inner solver.

		{\bf Funding.} The first author acknowledges financial support from the Austrian Science Fund (FWF) under \url{https://doi.org/10.55776/PAT2747625}. 
		
		{\bf Acknowledgment.} The authors would like to thank Christophe Tribes for helpful discussions concerning \texttt{NOMAD-4.6}, in particular its use, configuration, and interpretation of termination outcomes for linearly constrained optimization problems.
		
		\appendix
		\section{Appendix: Inner Solvers and Comparison Methods}
		\label{sec.inner-solvers}
		
		Subsections~\ref{sec:un+bd}--\ref{sec:LC} survey DFO methods for linearly constrained and unconstrained/bound-constrained problems. This appendix complements that discussion by identifying the specific routines used in the numerical study and summarizing their native algorithmic and constraint-handling capabilities.
		
		Table~\ref{a.sol} lists the routines used as inner solvers for the
		unconstrained transformed problem \gzit{e.unprob}, together with
		the proposed \texttt{VLDFO} framework in the last row.
		Some of these routines also support linear constraints and are
		used as direct comparison solvers for the original problem
		\gzit{e.lcopt}. The constraint types reported in the table describe
		the native capabilities of the solvers; when used as inner solvers,
		the routines are applied to the unconstrained reformulation.
		The table is not intended to provide an exhaustive survey of DFO
		solvers. Many other effective DFO methods are available and could
		also be used within the proposed framework. Since the present
		numerical study already includes a diverse collection of methods
		with different algorithmic characteristics, we restrict the table
		and comparisons to the solvers considered in the experiments.
		The omission of other methods should therefore be understood only
		as a limitation of the scope and size of the numerical comparison,
		and not as an assessment of their quality or applicability.
		
		The \texttt{fmincon-sqp} and \texttt{fmincon-active} entries select
		the SQP and active-set algorithms of the \texttt{MATLAB}
		Optimization Toolbox \cite{MatlabOTB}. In the inner-solver
		interface, their constraint arguments are empty and objective
		gradients are estimated by finite differences. Thus, they can be
		called using only function values, although the underlying
		algorithms remain gradient-based methods using numerical
		differentiation.

		\begingroup
		\setlength{\LTcapwidth}{\textwidth}
		\begin{longtable}{|lllllllll|}
			\hline
			solver & \rot{model-based} & \rot{line search} & \rot{direct search} &
			\rot{evolution strategy} & \rot{deterministic} & \rot{randomized} &
			\rot{constraint type} & Reference\\
			\hline
			\endfirsthead
			
			\multicolumn{9}{c}{\tablename~\thetable{} (continued)}\\
			\hline
			solver & \rot{model-based} & \rot{line search} & \rot{direct search} &
			\rot{evolution strategy} & \rot{deterministic} & \rot{randomized} &
			\rot{constraint type} & Reference\\
			\hline
			\endhead
			
			\hline
			\multicolumn{9}{r}{\small Continued on the next page}\\
			\endfoot
			
			\endlastfoot
			
			\texttt{BCDFO} & $+$ & $-$ & $-$ & $-$ & $+$ & $-$ & u/b &
			\sca{Gratton} et al. \cite{GraTT}\\
			\texttt{UOBYQA} & $+$ & $-$ & $-$ & $-$ & $+$ & $-$ & u &
			{\sc Powell} \cite{Pow02}\\
			\texttt{NEWUOA} & $+$ & $-$ & $-$ & $-$ & $+$ & $-$ & u &
			{\sc Powell} \cite{Pow08}\\
			\texttt{DFOTR} & $+$ & $-$ & $-$ & $-$ & $+$ & $-$ & u &
			{\sc Bandeira} et al. \cite{Bandeira2012DFOTR}\\
			{\tt LINCOA} & $+$ & $-$ & $-$ & $-$ & $+$ & $-$ & u/b/l &
			{\sc Powell} \cite{LINCOA}\\
			{\tt COBYLA} & $+$ & $-$ & $-$ & $-$ & $+$ & $-$ & u/b/l/n &
			{\sc Powell} \cite{COBYLA}\\
			\hline
			\hline
			\texttt{VRBBO} & $-$ & $+$ & $-$ & $-$ & $+$ & $+$ & u &
			{\sc Kimiaei \& Neumaier} \cite{KimN}\\
			\texttt{VRDFON} & $+$ & $+$ & $-$ & $-$ & $+$ & $+$ & u &
			{\sc Kimiaei} \cite{Kim}\\
			\texttt{SSDFO} & $-$ & $+$ & $-$ & $-$ & $+$ & $+$ & u &
			{\sc Kimiaei} et al. \cite{KimNF}\\
			\texttt{SDBOX} & $-$ & $+$ & $-$ & $-$ & $+$ & $-$ & u/b &
			{\sc Lucidi \& Sciandrone} \cite{LucS}\\
			\texttt{DFLBOX} & $-$ & $+$ & $-$ & $-$ & $+$ & $-$ & u/b &
			{\sc Liuzzi} et al. \cite{Liuzzi2011}\\
			\texttt{DFNDFL} & $-$ & $+$ & $-$ & $-$ & $+$ & $-$ & u/b &
			{\sc Giovannelli} et al. \cite{Giovannelli2022DFNDFL}\\
			\texttt{FMINUNC} & $-$ & $+$ & $-$ & $-$ & $+$ & $-$ & u &
			\texttt{MATLAB} Optimization Toolbox \cite{MatlabOTB}\\
			\texttt{fmincon-sqp} & $-$ & $+$ & $-$ & $-$ & $+$ & $-$ & u/b/l/n &
			\texttt{MATLAB} Optimization Toolbox \cite{MatlabOTB}\\
			\texttt{fmincon-active} & $-$ & $+$ & $-$ & $-$ & $+$ & $-$ & u/b/l/n &
			\texttt{MATLAB} Optimization Toolbox \cite{MatlabOTB}\\
			\hline
			\hline
			\texttt{DSPFD} & $-$ & $-$ & $+$ & $-$ & $-$ & $+$ & u/b/l &
			{\sc Gratton} et al. \cite{GraRVZ19}\\
			\texttt{NMSMAX} & $-$ & $-$ & $+$ & $-$ & $+$ & $-$ & u &
			{\sc Higham} \cite{Hig}\\
			\texttt{NELDER} & $-$ & $-$ & $+$ & $-$ & $+$ & $-$ & u &
			\sca{Nelder \& Mead} \cite{NelderMead1965}\\
			\texttt{ADSMAX} & $-$ & $+$ & $+$ & $-$ & $+$ & $-$ & u &
			{\sc Higham} \cite{Hig}\\
			\texttt{MCS} & $+$ & $-$ & $+$ & $-$ & $+$ & $-$ & u/b &
			\sca{Huyer \& Neumaier} \cite{HuyN0}\\
			\texttt{SID-PSM} & $+$ & $-$ & $+$ & $-$ & $+$ & $-$ & u/b/l/n &
			{\sc Cust\'odio} et al. \cite{CusV,CusRV}\\
			\texttt{NOMAD-3.9} & $+$ & $-$ & $+$ & $-$ & $+$ & $+$ & u/b/l/n &
			{\sc Le Digabel} \cite{Dig11}\\
			\texttt{NOMAD-4.6} & $+$ & $-$ & $+$ & $-$ & $+$ & $+$ & u/b/l/n &
			{\sc Audet} et al. \cite{Audet2022NOMAD4}\\
			\hline
			\hline
			\texttt{CMAES} & $-$ & $-$ & $-$ & $+$ & $-$ & $+$ & u/b &
			\sca{Auger \& Ha\-nsen} \cite{AugH}\\
			\texttt{fMAES} & $-$ & $-$ & $-$ & $+$ & $-$ & $+$ & u &
			{\sc Beyer} \cite{Beyer2020}\\
			\texttt{GCES} & $+$ & $-$ & $-$ & $+$ & $-$ & $+$ & u/b/l/n &
			{\sc Diouane} et al. \cite{DioGV0,DioGV}\\
			{\tt MADFO} & $-$ & $+$ & $-$ & $+$ & $+$ & $+$ & u &
			{\sc Kimiaei \& Neumaier} \cite{MADFO}\\
			\texttt{DMADFO} & $-$ & $+$ & $-$ & $+$ & $-$ & $+$ & u &
			\sca{Kimiaei \& Babaie--Kafaki} \cite{DMADFO}\\
			\texttt{MATRS} & $+$ & $+$ & $-$ & $+$ & $-$ & $+$ & u/b &
			\sca{Kimiaei \& Neumaier} \cite{MATRS}\\
			\hline
			\hline
			\texttt{VLDFO} & $+$ & $+$ & $-$ & $+$ & $+$ & $+$ & l &
			this paper\\
			\hline
			
			\caption{Solvers used as inner routines and, where applicable, as
				direct constrained comparison methods, together with the proposed
				\texttt{VLDFO} framework. The {\it constraint type} describes native
				solver capabilities: unconstrained (u), bound constrained (b),
				linearly constrained (l), and nonlinearly constrained (n). The
				remaining columns indicate the algorithmic components and whether
				the cited solver variants are deterministic, randomized, or both.}
			\label{a.sol}\\
		\end{longtable}
		\endgroup

		\section{Appendix: Smoothness Assumptions and Nonsmooth Applicability}
		\label{sec.nonsmooth-applicability}
		
		The transformed objective $e_\gamma$ is locally Lipschitz continuous
		and generally nonsmooth. It is therefore important to distinguish
		between a solver that can be applied numerically to $e_\gamma$ and
		one whose convergence analysis directly covers nonsmooth objectives.
		
		For several model-based and line-search routines, the available
		convergence analyses assume smoothness; examples include
		\texttt{VRBBO} \cite{KimN} and \texttt{SSDFO} \cite{KimNF}.
		Such methods can nevertheless be applied to $e_\gamma$ because their
		implementations use function values, but their smooth convergence
		results do not by themselves establish convergence for the
		nonsmooth transformed problem. The same distinction applies to
		classical interpolation-based trust-region methods and to the
		\texttt{MATLAB} routines that estimate derivatives by finite
		differences.
		
		Other methods were developed to accommodate nonsmooth black-box
		objectives. In particular, \texttt{DFNDFL}
		\cite{Giovannelli2022DFNDFL} incorporates continuous-search
		mechanisms designed for nonsmooth optimization. The \texttt{MADS}
		framework implemented in \texttt{NOMAD} \cite{Nomad} uses mesh
		refinement and polling directions whose refining subsequences
		satisfy density properties supporting nonsmooth stationarity
		results. References for \texttt{NOMAD 3.9} and
		\texttt{NOMAD 4.6} are \cite{Dig11} and
		\cite{Audet2022NOMAD4}, respectively.
		
		Methods such as \texttt{CMAES} \cite{AugH} and related evolution
		strategies rely on population sampling, selection, and distribution
		adaptation. These mechanisms do not require explicit gradients,
		subgradients, or bundle models.
		
		Accordingly, all routines in Table~\ref{a.sol} can be assessed
		numerically as inner solvers for the transformed problem, but their
		theoretical interpretation must distinguish between methods whose
		analysis assumes smoothness and methods whose theory explicitly
		allows nonsmooth objectives.

	\end{sloppypar}
	

	\end {document}